\documentclass[11pt]{article}

\usepackage[margin=1in]{geometry}
\usepackage{amsfonts,amssymb,amsmath,amsthm}
\usepackage{graphicx}
\usepackage{subcaption}
\usepackage{psfrag}
\usepackage{calc}
\usepackage{epic}
\usepackage[english]{babel}
\usepackage{hyperref}
\usepackage{float}
\usepackage{latexsym}
\usepackage{color}
\allowdisplaybreaks

\usepackage{hyperref}
\newcommand{\footremember}[2]{%
    \footnote{#2}
    \newcounter{#1}
    \setcounter{#1}{\value{footnote}}%
}
 
\title{Discrete-time TASEP with holdback
}
\author{%
  Seva Shneer\footremember{HW}{Heriot-Watt University, V.Shneer@hw.ac.uk}%
  \and Alexander Stolyar\footremember{UIUC}{University of Illinois at Urbana-Champaign, stolyar@illinois.edu}%
  }
\date{}
\newtheorem{theorem}{Theorem}
\newtheorem{lemma}[theorem]{Lemma}
\newtheorem{proposition}[theorem]{Proposition}

\newtheorem{conjecture}[theorem]{Conjecture}
\newtheorem{definition}{Definition}

\theoremstyle{remark}
\newtheorem{remark}[theorem]{Remark}

\newcommand{\be}{ \begin{equation}}
\newcommand{\ee}{\end{equation}}
\newcommand{\ben}{ \begin{equation*}}
\newcommand{\een}{\end{equation*}}

\newcommand{\beql}[1]{\begin{equation}\label{#1}}
\newcommand{\eeql}{\end{equation}}
\newcommand{\eqn}[1]{(\ref{#1})}

\def\E{{\mathbb E}}
\def\P{{\mathbb P}}
\def\I{{\mathbb I}}
\def\R{{\mathbb R}}
\def\bbZ{{\mathbb Z}}

\def\L{{\mathcal L}}

\begin{document}

\date{\today}

\maketitle

\begin{abstract} 
 We study the following interacting particle system. 
 There are $\rho n$ particles, $\rho < 1$, moving clockwise (``right''), in discrete time, on $n$ sites arranged in a circle. 
 Each site may contain at most one particle. At each time, a particle may move to the right-neighbor site according to the following rules.
 If its right-neighbor site is occupied by another particle, the particle does not move.
 If the particle has unoccupied sites (``holes'') as neighbors on both sides, it moves right with probability $1$. If the particle has a hole as the right-neighbor and an occupied site as the left-neighbor, it moves right with probability $0<p<1$. (We refer to the latter rule as a ``holdback'' property.) From the point of view of holes moving counter-clockwise, this is a zero-range process.
 
 The main question we address is: what is the system steady-state flux (or throughput) 
 when $n$ is large, as a function of density $\rho$? 
 The most interesting range of densities is $0\le \rho < 1/2$. We define the system {\em typical flux} as the limit in $n\to\infty$ of the steady-state flux in a system subject to additional random perturbations, when the perturbation rate vanishes. Our main results show that: (a) the typical flux is different from the formal flux, defined as the limit in $n\to\infty$ of the steady-state flux in the system without perturbations, and (b) there is a phase transition at density $h=p/(1+p)$. If $\rho<h$, the typical flux is equal to $\rho$, which coincides with the formal flux. If $\rho>h$, a {\em condensation} phenomenon occurs, namely the formation and persistence of large particle clusters; in particular, the typical flux in this case is $p(1-\rho) < h < \rho$, which differs from the formal flux when $h < \rho < 1/2$.
 
 Our results include both the steady-state analysis (which determines the typical flux) and the transient analysis. In particular, we 
 derive a version of the Ballot Theorem, and show that the key ``reason'' for large cluster formation for densities $\rho > h$ is described by this theorem. 
 
 \end{abstract}
 
 

\section{Introduction}


\subsection{The model and motivation.}

The basic model, which is the focus of this paper, is the following
interacting particle system.
There are $n$ sites (or nodes), arranged in a circle and numbered from $0$ to $n-1$ in the clockwise order. There is a constant number $\rho n$ 
of particles in the system, where $\rho \in (0,1]$ is the particle density. There is at most one particle at each site at any time.
The system evolves in discrete time $t=0,1,2,\ldots$. The state at a given time is described as a sequence of particles (occupied sites) and holes (empty sites). We often refer to the clockwise and counter-clockwise directions as ``right'' and ``left,'' respectively. 
The particles never move counter-clockwise. Their clockwise movement, at each time, is governed by the following rules:
\begin{itemize}

\item[(a)] if a particle has another particle as a right-neighbor, it does not move;
 
\item[(b)] if a particle has holes as neighbors on both sides, it moves to the right-neighbor site with probability $1$; 

\item[(c)] if a particle has a hole as right-neighbor and a particle as left-neighbor, it moves to the right-neighbor site with fixed
probability  $p \in (0,1]$.

\end{itemize}


This model may be considered as a version of the discrete-time Totally Asymmetric Simple Exclusion Process (TASEP),
with parallel updates, cf. \cite{Evans1997}. We refer to rule (c) as a ``holdback'' property and thus call the model discrete-time TASEP-H, where H stands for holdback. In the classical version of discrete-time TASEP a particle cannot move if its right-neighbor is another particle
and, otherwise, moves right with a certain fixed probability.
(For a general introduction into interacting particle systems see \cite{Liggett2005}.).


This interacting particle system is motivated, in particular, by a simple model of packet movement in a wireless communication system 
under a CSMA (Carrier-Sense Multiple Access) protocol. Sites correspond to network nodes and particles correspond to data packets. 
Discrete time corresponds to the sequence of fixed-length time slots, in which packet transmissions occur.
Each node can hold at most one packet.
Packets ``move'' along a sequence of nodes, in the ``right'' direction, by being transmitted from nodes to their ``right'' neighbors. 
Consider a node that in a given time slot has a packet. If both neighbors of the node do not have packets, the packet will be transmitted successfully, because the transmission does not experience an interference (from neighbors). If the right-neighbor node is occupied by another packet, the packet will be certainly blocked from moving right.
If the right-neighbor node is empty, but the left-neighbor node is occupied, the transmission protocol is such that the ``competition'' from the left-neighbor packet may prevent the packet from moving right, with some probability. (This motivating model 
is very basic. A more complicated -- and more realistic -- CSMA model, which exhibits a qualitatively similar behavior, is discussed in 
Appendix \ref{subsec:stork}.) Our TASEP-H system is also motivated by and related to a slow-to-start model \cite{CFP2007} of car traffic, where particles represent cars and the holdback property corresponds to the fact that cars need some time to accelerate after being stopped. 
The relation of TASEP-H to the slow-to-start model will be discussed in more detail in Section~\ref{sec-sts}.

Going back to the TASEP-H system definition,
note that from the point of view of holes moving counter-clockwise (left), the rules (a)-(c) are equivalent to the following ones:

\begin{itemize}

\item[(a')] if a hole has another hole immediately to the left, it does not move;
 
\item[(b')] if a hole has a particle immediately to the left, followed by another hole, it moves to the left-neighbor site;

\item[(c')] if a hole has two or more particles immediately to the left, it moves to the left-neighbor site with 
probability $p \in (0,1]$.

\end{itemize}

In the latter interpretation (as holes' movement couter-clockwise), TASEP-H system is a special case of the model considered 
in \cite{Evans1997} where the probability of a hole movement may be an arbitrary function of the number of holes immediately to the left of it.
It is shown (see \cite[Section 7]{Evans1997}) that the stationary distribution of the system may be presented in a product form.

In general, the class of TASEP models is very rich and has received a lot of attention in literature, because it has many applications, including statistical physics and transportation; see, e.g.\cite{Liggett2005,Blythe2007,Chowdhury2000,Helbing2001}.
In addition to \cite{Evans1997}, our model is related to some other TASEP models considered in literature. We mention $q$-TASEP  \cite{Borodin2014} for the continuous-time and \cite{Borodin2015} for the discrete-time version; both papers however consider movement of a finite number of particles on an infinite line. 

\subsection{Particle flux: ``formal'' Vs ``typical.''}

The particle {\em flux}  $\phi$ is defined as a steady-state average number of particle movements, per site per time unit. 
(This corresponds to the throughput in the context of wireless systems.)
In this paper we are primarily interested in the dependence of flux $\phi$ on the density $\rho$, when $n$ is large.

\begin{figure}
\centering
\includegraphics[width=0.7\linewidth]{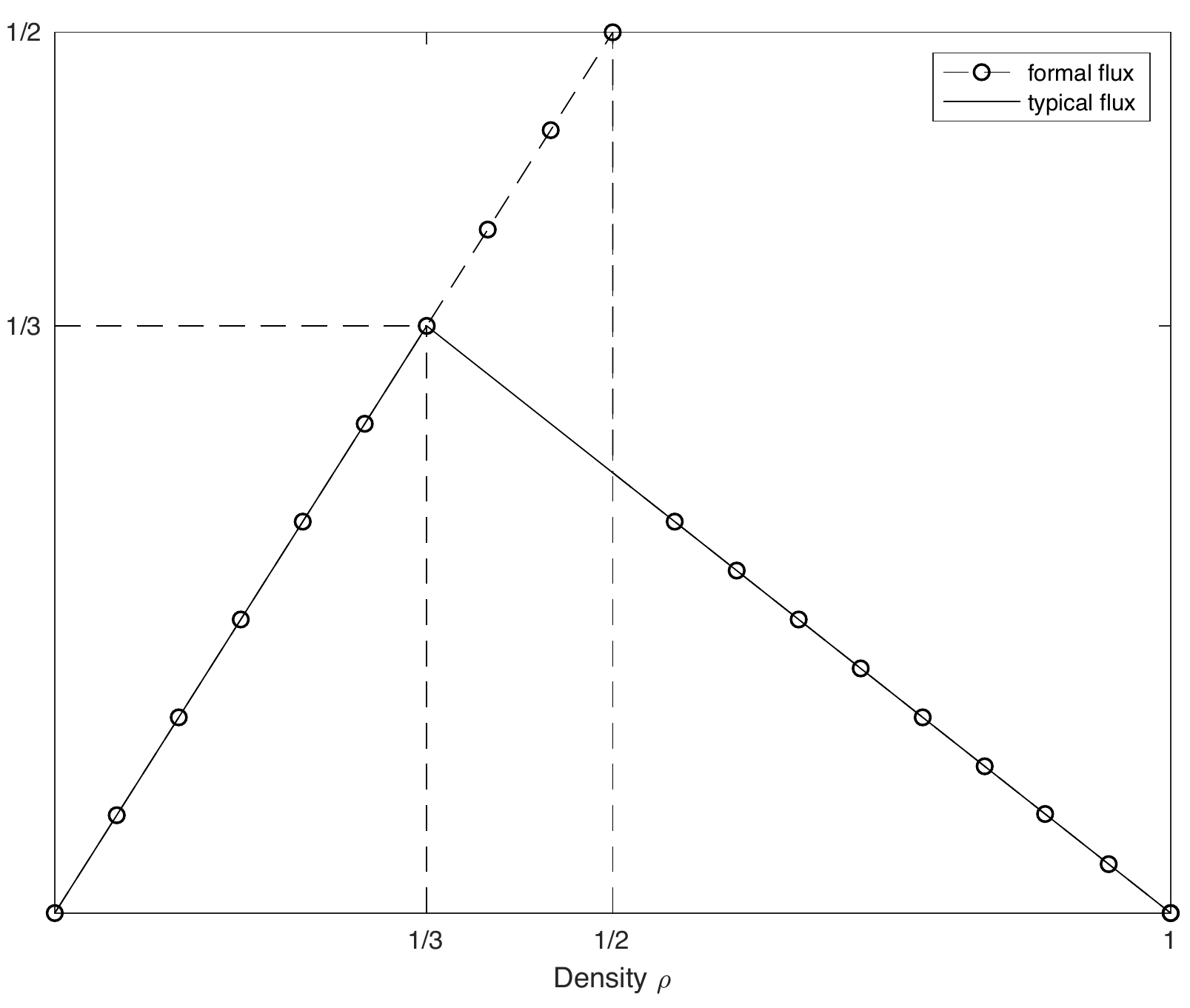}
\caption{Large $n$; $p=1/2$. Formal and typical flux versus $\rho$. For $\rho<1/3$ and $\rho>1/2$, the formal and typical fluxes coincide.}
\label{fig:flux_tasep_h_real}
\end{figure}

Trivially, for any $n$, when $\rho < 1/2$, the flux is simply equal to the density, $\phi=\rho$,
because the system eventually enters a {\em completely sparse} state, where all particles are ``free'' (never have neighbors),
and all particles will move at speed $1$ thereafter. 
Following the discussion which will be given shortly, 
it is also not hard to guess (and observe in simulations, and prove) that
when $\rho>1/2$, the steady-state of the system is very different -- the system spends most of the time in a {\em condensed} state, where a large particle cluster (a sequence of contiguously occupied sites) exists. Correspondingly, the ``formal'' flux of TASEP-H is as follows: $\phi=\rho$ for $\rho<1/2$, but it has has a discontinuity (``negative jump'') at $1/2$ and is continuous decreasing for $\rho>1/2$. Figure~\ref{fig:flux_tasep_h_real} shows this ``formal'' flux, when $n$ is large. 

However, an interesting phenomenon can be observed in simulations. If $h=p/(1+p) < \rho<1/2$ and $n$ is large, and we start the process in a state chosen uniformly at random, then the process enters a quasi-stationary regime, in which it stays in condensed states (with one or more clusters). This quasi-stationary regime persists for a very long time. The flux corresponding to this quasi-steady-state is strictly less than $\rho$.
If we consider as ``typical'' the system flux that is observed with very high probability for a very long time, then this ``typical'' flux of TASEP-H is shown in Figure~\ref{fig:flux_tasep_h_real} (for large $n$): it depends on $\rho$ continuously, coincides with ``formal'' flux for $\rho < h=p/(1+p)$ and $\rho>1/2$, but is different for $h < \rho < 1/2$.

In this paper we will formally define the ``typical'' flux as the limit in $n\to\infty$ of the flux of the TASEP-H system subject to additional state 
perturbations, as the perturbation rate (per particle per time unit) vanishes as $n\to\infty$. Our results will show, in particular, that the typical 
flux of TASEP-H is indeed as shown in Figure~\ref{fig:flux_tasep_h_real}.

\subsection{Condensation. Quasi-stationary regime.}

The term {\em condensation}  in interacting particle systems is inspired by the famous Bose condensation  \cite{Evans1996, Evans1997}.
(Regarding the terminology, in most of the literature, the condensation occurs at {\em low} particle density, as opposed to high.
This is consistent with TASEP-H condensation, if we view the process from the point of view of holes' movement.)

We now discuss  the condensation effect 
in more detail. Recall that the state of the system at any time is described as a  sequence of particles and holes. 
We will refer to a contiguous segment of (at least two) particles as a cluster. Only the particle at the right end of a cluster may move at a given time. We will also refer to a contiguous segment of sites (which may include all sites) as a sparse interval, if all particles in it are ``free,'' i.e. have only holes as neighbors. Obviously, each free particle at a given time moves right with probability $1$.

A system state such that all particles are free we will call (interchangeably)  as {\em completely sparse}, or {\em ideal}, or {\em absorbing}.
As discussed earlier, if $\rho<1/2$, then for any $n$ and any $p$, the system eventually reaches an absorbing state and therefore its flux $\phi=\rho$. 

However, for $h< \rho < 1/2$, it may take the system a very long time to reach such an absorbing state. Instead, the system may enter a {\em quasi-stationary} regime, which will persist for a long 
time. (It is reasonable to guess that this time is exponential in $n$, but we will not need this estimate in this paper.) 
We explain this -- heuristically at this point -- as follows. Suppose, the system starts in a state such that $\tau^* n$ particles, $0<\tau^* \le 1$, form a single cluster, while the the remaining sparse interval, consisting of $(1-\tau^*) n$ sites, has particles spread out with the density 
$$
h=p/(1+p).
$$ 
Note that in this case the overall density of particles is $\rho= \tau^* + (1-\tau^*)h > h$, and then 
$$
\tau^* = \rho(1+p) - p.
$$
Then, the system dynamics is as follows. 
The particles will leave the cluster from its ``right end'' at the rate $p$, with the average distance between the released particles being $1+1/p = 1/h$. Therefore, the cluster right end, as it moves left at the rate $p$, leaves particle density exactly $h$ in ``its wake'' on the right. Of course, the free particles (in the sparse interval) move right at speed $1$. It is not hard to see that, as free particles ``hit'' (and thus join) the cluster on its left, the left end moves left at the rate $p$. (This follows because the density in the sparse interval is $h$.) To summarize, the cluster moves left at speed $p$, with its length $\tau^* n$ staying (approximately) constant, and with the density ``everywhere'' in the sparse interval staying equal to $h$. Therefore, the system state remains (almost) invariant, up to the shift of the cluster. We will call such a state a {\em main equilibrium state} (MES). (This definition will be made formal later in the paper.) Note that a MES exists and is unique (up to space shift) if and only if $\rho > h$. As long as a MES remains the system state,
the flux is $\phi^* \doteq (1-\tau^*) h = p(1-\rho) < h < 1/2$. The existence of a quasi-stationary regime, represented in our model by MES, is often referred to as {\em metastability}.

Further notice that, if the system starts in a MES, it will take a long time (which, again, can be guessed to be exponential in $n$) for the cluster to ``dissolve.''
Therefore, if the density $\rho \in (h, 1/2)$, the system may enter a MES and stay in it
 -- as opposed to an absorbing state -- for a very long time, during which it will have the flux  $\phi^* < h < \rho$, as opposed to the ``formal'' flux $\rho$. 

\begin{figure}
\centering
\includegraphics[width=0.7\linewidth]{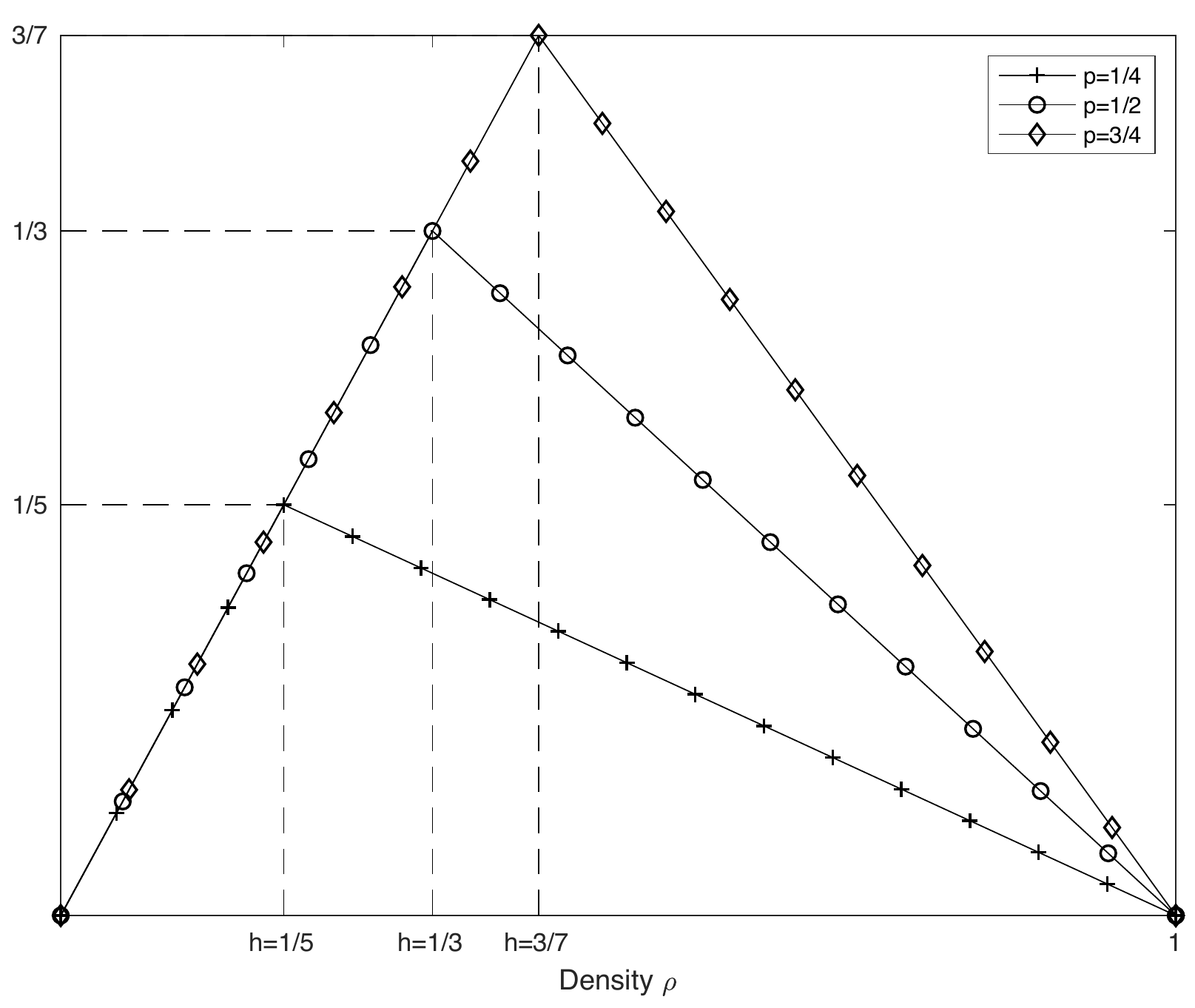}
\caption{Large $n$; $p=1/4, 1/2, 3/4$. Typical flux versus density $\rho$}
\label{fig:flux_tasep_h_typical_multiple}
\end{figure}

Note that MES exists for any $\rho > h$, not just $h < \rho < 1/2$. 
 The graph of $\phi^*$ (when the system stays in a MES) versus $\rho$, for $\rho> h$,  can be seen in Figure \ref{fig:flux_tasep_h_typical_multiple}, for the cases $p=0.1$, $p=0.5$ and $p=0.9$. Note that $\phi^*$  is linearly 
 decreasing from $h$ to $0$ in the interval $[h,1]$.
 The case when $\rho \ge 1/2$ is a simplified version of the case $h < \rho < 1/2$. For this reason, in this paper we 
 restrict our attention to the most interesting case $\rho < 1/2$.

\subsection{Main results and contributions}

To summarize the previous subsection, even when $\rho< 1/2$, in addition to a formal stationary regime, when the system stays in an absorbing state and the flux is $\rho$, there may exist a quasi-stationary regime (represented by a MES), in which the system flux is $\phi^*$.
This raises the following question. What should be considered a ``typical'' flux of TASEP-H?
Is it $\rho$, or $\phi^*$, or maybe something else? 
In this paper we introduce the notion of typical flux, which is, informally speaking, the value of the flux in the system subject to infrequent random perturbations 
of the system state, in addition to the system ``normal'' random evolution. 
This is motivated by the fact that, even if the basic system behavior is described by TASEP-H -- whether it is movement of packets in a wireless system or movement of cars on a road -- practical systems are subject to occasional ``imperfections,'' or perturbations.
For example, new packets [resp., cars] may occasionally enter or leave the system [resp., road]; packets may occasionally fail to transmit successfully even in the absence of interference; cars may break down; etc. If occasional perturbations result in the flux that is different from the formal flux, then the former rather than latter can be considered typical.
More specifically, in this paper, for a given density $\rho$, we will define the typical flux $\phi(\rho)$ as the limit in $n\to\infty$ of the flux of
TASEP-H {\em with perturbations}, as the rate of perturbations (per particle per time unit) vanishes as $n\to\infty$.
The main goal of this paper is to demonstrate that even very low perturbation rate, $O(1/n^2)$, leads to the typical flux given by Figure \ref{fig:flux_tasep_h_real}.

Perturbations in the basic TASEP-H may be defined in different ways. For example, we may replace the probability $1$ with which 
a free particle moves right, with a close to one probability $\pi<1$. Such a system,
which we will refer to as {\em zero-range} model, 
is still within the framework of \cite{Evans1997}, its stationary distribution has a product form, and the limiting flux $\phi_{\pi}$, as $n\to\infty$, can be found. 
We in fact carry out this analysis in
Section~\ref{sec-zrp}, and demonstrate that $\lim_{\pi\uparrow 1} \phi_{\pi}$ is exactly as depicted on Figure \ref{fig:flux_tasep_h_typical_multiple}.
Note that here we take the limit $n\to\infty$ first, and limit $\pi\uparrow 1$ second. This means, in particular, that as $n\to\infty$, the perturbation rate remains ``high,'' namely $O(1)$ per particle per time unit. 

However, the above analysis, based on the steady-state product-form for the zero-range model has a number of limitations, including the following.

\begin{itemize}

\item It is interesting to look at how low a perturbation rate can be so that 
the flux remains as in Figure \ref{fig:flux_tasep_h_typical_multiple}. This question could be answered by considering a zero-range model with $\pi=\pi(n)\uparrow 1$. 
Unfortunately, analysis of the flux in this latter case (i.e., taking limit $n\to\infty$ with $\pi=\pi(n)\uparrow 1$), appears to be difficult.

\item Perturbations induced by $\pi < 1$ are not the only type of perturbations of interest. And for other types of perturbations, as those we introduce shortly, the product-form for the steady-state may no longer be available.

\item An analysis of steady-state alone does not shed light into how condensation occurs, i.e. into the questions of the following type. What is the mechanism and dynamics of large cluster formation?
How fast do the clusters form? An so on.

\end{itemize}

To address these issues, most of the analysis in this paper
concerns a different type of perturbations (introduced shortly) and ``low'' perturbation rates (per particle per time unit), vanishing as $n\to\infty$. For these perturbations the product form of the steady-state is {\em not} available.
In particular, we do a transient analysis which sheds light on why, for the densities $\rho > h$, even a small cluster formed by a perturbation has a positive, uniformly bounded away from $0$, probability to grow large, bring the system to a MES, which will then persist for a long time.
The key ``reason'' for large cluster formation for densities $\rho > h$ is described by a version of the Ballot Theorem, which we derive
in Proposition~\ref{th-ballot}.

The TASEP-H with perturbations, which is the main focus of this paper, is as follows.
We define two versions of the perturbation mechanism, both preventing the system from being absorbed in an ideal (absorbing) state. For both of them, a perturbation is defined as follows: a particle, chosen uniformly at random, is relocated to one of the holes, also chosen uniformly at random. The two versions of TASEP-H with perturbations are:

\begin{itemize}

\item {\it A-perturbations (perturbations of absorbing states).} At every time, in which the system enters an ideal (absorbing) state, a perturbation happens with probability $\lambda/n$ with some $\lambda > 0$, independent of the process history up to that time.

\item {\it I-perturbations (independent perturbations).} At every time a perturbation happens with probability $\lambda/n$ with some $\lambda > 0$, independently of the process history up to that time.

\end{itemize}

We analyze the behavior of the system as $n\to\infty$, with $\rho$ staying constant. In particular, we study the process under fluid (mean-field) space/time scaling. In addition to the steady-state results -- in fact, as a tool for obtaining them -- we derive results on the process transient 
behavior under this scaling. These transient results may be of independent interest.

Our {\bf main results} 
are as follows.

\begin{itemize}

\item
For {\bf TASEP-H with A-perturbations} we prove that the limit (in $n\to\infty$) of stationary distributions is such that: 

\begin{itemize}

\item If $\rho < h$, a cluster (if any) contains zero fraction of sites; consequently, the typical flux $\phi(\rho)=\rho$; see Theorem~\ref{th-A-low}.
\item If $\rho > h$, the system state is the MES; consequently, the typical flux $\phi(\rho)=\phi^*=p(1-\rho)$; see Theorem~\ref{th-A-high}.
\end{itemize}

\item
For {\bf TASEP-H with I-perturbations} we prove that the limit (in $n\to\infty$) of stationary distributions is such that:
\begin{itemize}
\item If $\rho < h$, all clusters (if any) contain zero fraction of sites; consequently, the typical flux $\phi(\rho)=\rho$; see Theorem~\ref{th-I-low}.
\item If $\rho > h$, the typical flux $\phi(\rho) < \rho$; see Theorem~\ref{th-perturb}.
\end{itemize}
We conjecture that, if $\rho > h$, the typical flux $\phi(\rho)=\phi^*$; Conjecture~\ref{conj-high-dens}.

\item
For the {\bf zero-range model} we prove that if we first take the limit of stationary distributions as $n\to\infty$, and then the limit in $\pi\uparrow 1$, the limit is exactly equal to the typical flux $\phi(\rho)$ of TASEP-H with A-perturbations; see Theorem~\ref{th-zrp} in Section~\ref{sec-zrp}.

\item
We derive a version of the {\bf Ballot Theorem} (Proposition~\ref{th-ballot}) which serves as a key tool for establishing condensation at high densities $\rho > h$.

\end{itemize}

\subsection{Paper organization}

The rest of the paper is organized as follows. Some basic notation, used throughout the paper, is given in Section~\ref{sec-notation}.
Section~\ref{sec-model} defines the basic TASEP-H model (without perturbations), the model with perturbations, and a system flux; it also defines the asymptotic regime, the corresponding limiting flux (typical flux), the process under fluid (mean-field) scaling, and the probability space construction.
The analysis of TASEP-H with A-perturbations is given in Section~\ref{sec-a-perturb}: we define the process fluid limits and establish their properties (Section~\ref{sec-a-perturb-fl}); we state a version of the Ballot Theorem (Proposition~\ref{th-ballot} in Section~\ref{sec-ballot-prop});
the main results  for TASEP-H with A-perturbations (Theorems~\ref{th-A-low} and \ref{th-A-high}) are stated and proved in 
Section~\ref{sec-A-pert-main-res}. Section~\ref{sec-i-perturb} is devoted to the analysis of TASEP-H with I-perturbations: 
some preliminary facts are given in Section~\ref{sec-i-prelim}; fluid limits are defined and studied in Section~\ref{sec-i-perturb-fl};
Sections~\ref{sec-i-low-result} and \ref{sec-i-pert} contain our main results and the conjecture (Theorem~\ref{th-I-low}, Theorem~\ref{th-perturb},
Conjecture~\ref{conj-high-dens}) for TASEP-H with I-perturbations. In Section~\ref{sec-zrp} we formally define the zero-range model and 
derive (Theorem~\ref{th-zrp}) the limit of its flux (as $n\to\infty$ and then $\pi \to 1$). In Section~\ref{sec-discuss-conj} we present and discuss further conjectures for both the
TASEP-H with perturbations and zero-range models. Appendix~\ref{app-ballot} contains the proof of Proposition~\ref{th-ballot};
Appendix~\ref{subsec:stork} discusses
a CSMA model further (indirectly) motivating the holdback property of TASEP-H.

\subsection{Basic notation}
\label{sec-notation}

We denote by $\R$ and $\bbZ$ the sets of real numbers and integers, respectively, and by
 $\L$ the Lebesgue measure on $\R$. Abbreviation {a.e.} means {\em almost everywhere w.r.t. Lebesgue measure}.
 Notation $(\partial^-/\partial x) f(x,t)$ means left partial derivative in $x$.
 The minimum and maximum of two numbers  are denoted 
 $a \wedge b \doteq \min(a,b)$ and $a \vee b \doteq \max(a,b)$, respectively. $\I(A)$ 
 is the indicator of event or condition $A$. RHS and LHS mean right-hand side and left-hand side, respectively.
 
 Abbreviation {\em w.p.1} means {\em with probability $1$}. Probability distributions are defined on the spaces and 
 corresponding $\sigma$-algebras that are clear from the context.
 We denote by $\Rightarrow$ the convergence of random elements in distribution.
 For a random process $Y(t), t\ge 0$, we denote by $Y(\infty)$ its (random) state in a stationary regime. (In other words, the distribution of $Y(\infty)$ is a stationary distribution of $Y(t)$.) RCLL means {\em right-continuous with left limits}.
We will say that a random variable $X$ has distribution $GEOM(\ell,p)$ [or simply write $X \sim GEOM(\ell,p)$], for integer $\ell$ and real $0\le p \le 1$,
 if $\P\{X=i\}= p(1-p)^{i-\ell}$, $i\ge \ell$; $\E X = \ell + (1-p)/p$.

\section{Model} 
\label{sec-model}

\subsection{Basic model}

Consider $n$ sites arranged in a circle.  The sites are numbered from $0$ to $n-1$ in clockwise order. 
There is a constant  number $\rho n$, $0 < \rho \le 1$, of particles in the system, with either $0$ or $1$ particles located at each 
site at any time; correspondingly, a site may be empty ($0$ particles) or occupied ($1$ particle). 
We often refer to empty sites as holes, and to occupied sites as particles.
The clockwise [resp., counterclockwise] direction we often refer to as ``right'' [resp. ``left''] direction. 

The system evolves in discrete time $t=0,1,2,\ldots$. A time index $t$ is sometimes referred to as time slot $t$.
The system state at each time is its sites' occupancy configuration, i.e. the sequence of particles and holes. Given the state at time $t$, the (random) state at time $t+1$ is determined by applying the following rules to each particle independently:\\
 (a) if the right-neighbor of the particle is another particle, it does not move; \\
 (b) if both the right-neighbor and left-neighbor of the particle are holes, the particle moves to the right-neighbor site with probability $1$; \\
(c) if the right-neighbor of the particle is a hole and the left-neighbor is another particle, the particle moves to the right-neighbor site with 
probability $p \in (0,1]$.\\
We refer to the rule (c) as a ``holdback'' property.

To make some definitions that follow later unambiguous, we adopt the following convention about the exact timing of the particles' movement. Namely, the movement of particles that changes the system state from that at $t-1$ to that at $t$ is attributed to time $t$.

We will use the following terminology throughout. We call a particle {\em free} if it has holes on both sides. A contiguous set of particles is a {\em cluster} if it has holes as neighbors on both sides.  A contiguous set of sites is a {\em sparse interval}, if it contains only free particles and holes (i.e., does not overlap with any cluster). A system state is called {\em completely sparse}, or {\em ideal}, if it contains no clusters.

Note that any trajectory of the basic model 
is such that the number of clusters cannot increase. Indeed, let us assume that as particles ``leave'' a cluster from the right and ``join'' it on the left, this cluster retains its ``identity.'' Then, each initial cluster may ``move'' to the left, may grow or decrease in length, or may eventually disappear. But no new cluster can ever be created.

\subsection{TASEP-H with perturbations}

We will now introduce the TASEP-H model with perturbations, in fact two different versions of it. For a given system state its (random) {\em perturbation} is defined as follows: we pick a particle uniformly at random, remove it, and then place it into one of the holes picked uniformly at random. We adopt a convention that the site from which the particle was removed is immediately considered a hole, so that the particle may go back to it. (This convention is not essential.)

We also adopt the following convention about the exact timing of a perturbation with respect to time slots.
Applying a perturbation at time $t$ means that it is applied to the state the system enters 
after the ``normal'' particles' movement at time $t$ (i.e. the movement ``between'' $t-1$ and $t$), and it is done before the next time $t+1$. 
In other words, the perturbation does not take an extra time slot, 
and the final system state at time $t$ is the state after the normal particles' movement at $t$ and after a perturbation (if any) is applied.

{\em Model with I-perturbations.} In this version, at each time $t$, one perturbation is applied with probability $\lambda/n$, $\lambda>0$, independently of the process history up to $t$. (Term I-perturbations is because the perturbations are independent.)

{\em Model with A-perturbations.} In this version, one perturbation is applied with probability $\lambda/n$, $\lambda>0$, 
at each time $t$ in which the system enters to -- or stays at -- an ideal state, {\em after the normal particles' movement.} (Term A-perturbations is because the perturbations explicitly prevent absorption.)

Recall that in the basic model (without perturbations) no new cluster can ever be created.
Unlike in the basic model, a model with perturbations is such that new clusters may be created (only after a perturbation). 
Note that, for the A-perturbation model, this necessarily means that, after the system ``hits'' and ideal state, it can have at most one cluster from that time on. In the I-perturbation model, a cluster may be split into two (only after a perturbation). 

For future reference, we summarize these simple observations as the following

\begin{lemma}
\label{lem-simple} Any process trajectory has the following properties, depending on the model type. 

(i) The basic model (without perturbations). The number of clusters cannot increase. Each initial cluster may move to the left, may grow or decrease in length, may eventually disappear. But no new cluster can ever be created.

(ii) A-perturbations model. The number of clusters cannot increase until a trajectory ``hits'' an ideal state. There can be at most one cluster from that time on, which may be created after a perturbation.

(iii) I-perturbations model. The number of clusters cannot increase between perturbations. Upon a perturbation, a new cluster may form and/or a cluster may split into two.
\end{lemma}

The process without perturbations is trivial in that an ideal state is eventually reached.
For both the A-perturbation and I-perturbation models,
the corresponding process is a discrete-time finite irreducible aperiodic Markov chain. 
For A-perturbations, the state space is irreducible if it is a priori restricted to those states reachable from the state where all particles are within one cluster (all other states are transient).

\subsection{Flux of TASEP-H with perturbations.} 

We now define the flux of TASEP-H with perturbations. Denote the state of the system at time $t$ by $\overline{Z}(t) = (Z_1(t), \ldots, Z_n(t))$, where $Z_i(t)=1$ if site $i$ is a particle and $Z_i(t)=0$ if it is a hole. 
Then, the instantaneous average flux of the system is defined as $\phi(\rho, n; \overline{Z}(t))=\psi(\rho, n; \overline{Z}(t))/n$, 
where $\psi(\rho, n; \overline{Z}(t))$ is the expected total distance that will be traveled by all particles at time $t+1$.
(This includes particles being relocated due to perturbations -- the precise convention will be given below.)
The {\em flux} of TASEP-H with perturbations, with parameters $\rho$ and $n$, is defined as
$$
\phi(\rho, n) \doteq \E \phi(\rho, n; \overline{Z}(\infty)),
$$ 
where $Z(\infty)$ is the (random) value of $Z(t)$ in a stationary regime. (It is easy to see that $\phi(\rho, n)$ for a given system 
does not depend on a steady state chosen, even if steady state is non-unique.)

\subsection{Asymptotic regime. Limiting flux (typical flux)} 

The asymptotic regime that we will consider is such that $\rho$ and $\lambda$ remain constant, while $n\to\infty$. To avoid cumbersome notation, let us assume that $\rho n$ is integer. This is the total number of particles. (If $\rho n$ is non-integer, we could assume that
the number of particles is, say, $\lceil \rho n \rceil$.)

The main focus of the paper is on identifying the {\em limiting flux} (of TASEP-H with perturbations):
$$
\phi(\rho) = \lim_{n \to \infty} \phi(\rho,n).
$$
The value of $\phi(\rho)$ will also be referred to as the {\em typical flux} (of TASEP-H). 

{\em Formal flux} (of TASEP-H) is the limit $\phi(\rho)$ for TASEP-H without perturbations.

\subsection{Process definition. Fluid (mean-field) scaling} 

For a system with parameter $n$ (the circle length) the process state at time $0$ is described by the function
$$
F^n(x,0), ~~ x=0,1,\ldots, n-1,
$$
where $F^n(x,0)$ is the total number of particles in sites $0, \ldots, x$ at time $0$. 

The (random) movement of particles in the system is described by the flux-function $\Phi^n(x,t)$, defined for 
$x=0,1,\ldots, n-1$ and time  $t=0,1,2, ...$ as follows: $\Phi^n(x,t)$ is the total number of particles that moved (right) from site $x$ at times $1, \ldots, t$, with $\Phi^n(x,0)=0$. By convention, this quantity includes the number of particles that crossed site 
$x$ due to perturbations; namely, the convention is that a particle being relocated from site $x_1$ to site $x_2$ always moves ``right,'' i.e. $x_1 \le x_2$, and this particle simultaneously ``leaves'' sites $x_1, x_1+1, \ldots, x_2-1$ (or none, if $x_1=x_2$).\footnote{Note that the contribution of perturbations into the steady-state flux is upper bounded by $n(\lambda/n)/n$, 
and therefore vanishes as $n\to\infty$.}

We define
$$
F^n(x,t) = F^n(x,0) - \Phi^n(x,t), ~~ x=0,1,\ldots, n-1, ~~ t=0,1,2, ...
$$
It is convenient to extend the definition of $F^n(x,t)$ to all integer $x$ by setting $F^n(n,t)=F^n(0,t)+\rho n$
and assuming that $F^n(x,t)$ has periodic increments in $x$ with period $n$.
Clearly, the (random) function $F^n(x,t)$ completely describes the process evolution. It is non-decreasing in $x$ and non-increasing in $t$. (In terms of describing the system states and evolution, only increments of $F^n(x,t)$ matter. So, 
for any fixed integer $C$, $F^n(\cdot,t)+C$ describes exactly the same system state at $t$ as $F^n(\cdot,t)$,
and $F^n(\cdot,\cdot)+C$ describes exactly the same system trajectory as $F^n(\cdot,\cdot)$. In particular,
if $F^n(\cdot,t_1)-F^n(\cdot,t_2)=C, ~\forall x$, the system states at times $t_1$ and $t_2$ are equal.)

It is also convenient to extend the definition of $F^n(x,t)$ to all real $x \in \R$ and $t\ge 0$ by adopting convention
$$
F^n(x,t) = F^n(\lfloor x \rfloor, \lfloor t \rfloor).
$$
Furthermore, we will identify any location $x\in \R$ 
with $x ~(\mbox{mod} ~n) \in [0,n)$.
So, for example, if $n=50$, $F^{50}(4,t) - F^{50}(4-7,t)$ is the total number of particles, at time $t$, in the $7$ 
consecutive sites, starting from $4$ and going left: $4,3,2,1,0,-1,-2$ or, equivalently, $4,3,2,1,0,49,48$. 

We assume that the system processes, for a given $n$, are constructed via driving (control) sequences of random variables, described as follows.

We have a countable set, indexed by $j=1,2,\ldots$, of the sequences of i.i.d. random variables, $\xi^j_i, ~i=1,2,\ldots,$ with distribution 
$GEOM(1,p)$, mean $1/p$.
For each cluster that exists initially, and for each new cluster that forms as the process evolves, a sequence 
$\{\xi^j_i\}$ with its ``own'' $j$ is assigned. So, $j$ can be thought of as a cluster index. The indices $j$ are assigned to the clusters sequentially, ``as needed.'' If a cluster $j$ breaks into two (which is possible due to and only due to a perturbation), then, by convention, the right one retains the index $j$, and the left one gets a new index. If two clusters $j'$ and $j$ merge into one larger cluster (which is possible due to and only due to a perturbation), 
with cluster $j$ becoming the right part of the new cluster, 
then, by convention, the new cluster retains the index $j$, and the index $j'$ is ``eliminated.''
A r.v. $\xi^j_i$ determines the random time it takes a particle to ``break away'' from a cluster, after this particle becomes the right edge of the cluster. (For example, {\em assuming} the particle just ahead of ours, that broke away from the same cluster earlier, remains free, then $\xi_i$ gives the number of holes between our particle and the particle ahead, at the time when our particle breaks away.)
The random variables $\xi^j_i$, for a given cluster $j$, are ``taken'' in sequence, as necessary, every time a particle becomes a cluster-right-edge particle. 

Perturbations are controlled by an i.i.d. sequence of pairs $(\eta_{1,i},\eta_{2,i}),~i=1,2,\ldots,$ with independent components, each having uniform distribution in $[0,1)$. Variable $\eta_{1,i}$ is ``responsible'' for picking a particle, and $\eta_{2,i}$ is ``responsible'' for picking a vacant site where it relocates.  Specifically, $\eta_{1,i}$ is ``used'' as follows. 
We label the particles by $1,2,\ldots, \rho n$, in the order of their locations in the (scaled) interval $[0,1)$. 
Then, the particle with label $\lceil \eta_{1,i} \rho n \rceil$ is picked for a perturbation. A variable $\eta_{2,i}$ is used analogously for picking a vacant site where a particle relocates. The timing of perturbations is controlled by the i.i.d. sequence $\eta_{3,i},~i=1,2,\ldots$, with
$\eta_{3,i} \sim GEOM(0,\lambda/n)$, $\E \eta_{3,i} = (1-\lambda/n)/(\lambda/n)$. Namely, $\eta_{3,i}$ determines the random time from the start of the $i$-th perturbation ``clock,'' until the actual $i$-th perturbation.  
For I-perturbations, the clock starts in the time slot immediately after the previous, $(i-1)$-th, perturbation;
for A-perturbations, the clock starts in the first time slot after the previous perturbation, where the system enters an ideal state.
Triples $(\eta_{1,i},\eta_{2,i}, \eta_{3,i})$ are ``used'' in sequence, as necessary, when perturbations need to be applied.

For each $n$, we consider a fluid (mean-field) scaled process, where we compress space, time
and the number of particles by factor $n$:
$$
f^n(x,t) = \frac{1}{n}F^n(nx,nt). 
$$
For a fixed (scaled) time $t$, $f^n(x,t)$  -- the (scaled) system state -- has periodic increments with period $1$. 
From now on, we always refer to scaled space and time, unless explicitly stated otherwise. 

In this paper we study the asymptotic behavior of the fluid-scaled process $f^n(x,t)$ as $n\to\infty$. 

In the rest of the paper, we will often construct the processes for all $n$ on a common probability space. The construction is as follows. We assume that the  i.i.d. sequences, $\xi^j_i, ~i=1,2,\ldots,$, $j=1,2,\ldots$, are common for all $n$. They satisfy the following functional strong law of large numbers (FSLLN) condition: w.p.1,
\beql{eq-xi-fslln}
\lim_{n\to\infty} \frac{1}{n} \sum_{i=1}^{\lfloor nt \rfloor} \xi^j_i = \frac{1}{p} t, ~~t\ge0 ~(u.o.c.), ~~~\forall j.
\eeql
The i.i.d. sequences
$\eta_{1,i}$ and $\eta_{2,i}$, $i=1,2,\ldots,$, will also be common for all $n$. The i.i.d. sequence $\eta_{3,i}^n$, $i=1,2,\ldots,$ will depend on $n$. Using Skorohod representation, we can and will assume that these sequences are such that, w.p.1, 
\beql{eq-eta-conv}
\lim_{n\to\infty} n^{-1} \eta_{3,i}^n = \hat \eta_{3,i}, ~~~\forall i,
\eeql
where $\hat \eta_{3,i}$ are independent (across $i$) exponentially distributed random variables with mean $1/\lambda$. Of course, we also have, by the strong law of large numbers, that w.p.1,
\beql{eq-eta-slln}
\lim_{m\to\infty} \frac{1}{m} \sum_{i=1}^m \hat \eta_{3,i} = \frac{1}{\lambda}.
\eeql

\section{A-perturbations model}
\label{sec-a-perturb}

\subsection{A-perturbations model: Fluid limits}
\label{sec-a-perturb-fl}

For any $n$, $f^n(\cdot,t)$ is a Markov process, and the main goal of this paper is to study its asymptotic behavior, in particular the asymptotics of its stationary distribution. However, to achieve this goal, we will need to analyze an extended version of this process, which includes
additional process variables, explicitly 
describing clusters, sparse intervals, and the times when clusters exist. 
Recall that, for the A-perturbations, without loss of generality, we restrict the state space to those states reachable from a state where all particles are within one cluster. Therefore, for A-perturbations model, there is at most one cluster at any time, so the following additional descriptors will suffice.

Denote by $\ell^n(t)$ and $r^n(t)$ the (scaled) locations of the particles at the left and right edge of the cluster at time $t$, if any, and by $\tau^n(t)=r^n(t)-\ell^n(t)$ its length.
(Note that the scaled number of particles within the cluster is  $[n\tau^n(t)+1]/n = \tau^n(t) + 1/n$, not $\tau^n(t)$.)
By convention, if cluster does not exist at time $0$, we assume that $\ell^n(t)=r^n(t)=\tau^n(t)=0$ until the time when a cluster forms.
Also by convention, after a cluster dissolves, $\ell^n(t)=r^n(t)$ remain frozen at the value of $\ell^n$ at the last time when the cluster existed.
We also adopt the convention that when $\ell^n$ changes as a result of a cluster formation, the new value of $\ell^n$ is chosen 
to be within $[0,1)$, and $r^n$ is chosen accordingly, to be to the right of (or at) $\ell^n$.

At any time $t$, the system state is given by
$$
[f^n(\cdot,t),\ell^n(t), r^n(t), \tau^n(t)].
$$
(Note that $[f^n(\cdot,t),\ell^n(t), r^n(t), \tau^n(t)]$ contains no more information than $f^n(\cdot,t)$, in the sense that the distribution of
$f^n(\cdot,s), s \ge t$ depends only on $f^n(\cdot,t)$. Once again, we use the additional state descriptors for the purposes of analysis.)
The metrics on the state space components are as follows. For the $f$-component, it is defined by the max-norm:
$$
\|f_1(\cdot)-f_2(\cdot)\| = \max_x |f_1(x)-f_2(x)|.
$$
For the $\tau$-component, it is given simply by the distance $|\tau_1 -\tau_2|$. 
For the $\ell$-component, it is the distance between
$\lfloor \ell_1 \rfloor$ and $\lfloor \ell_2 \rfloor$, considered as points on the unit length circle; formally, it is
$$
\min_{k\in \mathbb Z} |\ell_1 -\ell_2 +k|.
$$
The metric for the $r$-component space is the same as for the $\ell$-component.
Finally, the metric on the entire system state space is the sum-metric of its components.

For each $n$, consider a finite or infinite time interval $[0,T^n]$, with $T^n$ possibly depending on $n$. 
(If $T^n$ does depend on $n$, it can be thought of as a realization of a stopping time.)
Trajectories
$$
[f^n(\cdot,t),\ell^n(t), r^n(t), \tau^n(t)], ~t\in [0,T^n],
$$ 
are considered as elements of the Skorohod space. For the purposes of having the Skorohod metric well-defined, we adopt the convention that, when $T^n < \infty$, the state remains constant in $[T^n,\infty)$.

Further, let $\alpha^n_0\ge 0$  be the first time when a cluster exists. (If it exists at time $0$, then $\alpha^n_0=0$. If it never exists 
$\alpha^n_0=\infty$.) Let $\beta^n_0 \ge \alpha^n_0$ be the first time after $\alpha^n_0$ when the system is in an ideal state. 
($\beta^n_0=\infty$ if the system never enters an ideal state after $\alpha^n_0$.) In other words, $[\alpha^n_0,\beta^n_0)$ is the time interval when the ``first'' (in time) cluster exists. Similarly, let $\alpha^n_1\ge \beta^n_0$ be the first time after $\beta^n_0$ when the ``second'' cluster appears, 
and $\beta^n_1 \ge \alpha^n_1$ be the first time this cluster disappears. And so on, we define pairs $(\alpha^n_i,\beta^n_i)$
marking the beginning and end of the $i$-th cluster. Note that, $\alpha^n_i =\infty$ implies that all consecutive times $\alpha^n_j$ and $\beta^n_j$ are also infinite. We also adopt a convention that, if we consider the process 
on a finite time interval $[0,T^n]$, then any time $\alpha^n_i$ or $\beta^n_i$ which is outside $[0,T^n]$ is set to be infinite.
Finally, note that all $\alpha^n_i$, except maybe $\alpha^n_0$ when it is $0$, are the times when perturbations occur. (However,
not every perturbation time is one of the $\alpha^n_i$.)

The following will be called an extended realization (trajectory) of the process on a (finite or infinite) time interval $[0,T^n]$:
$$
\{ [f^n(\cdot,t),\ell^n(t), r^n(t), \tau^n(t)], ~t\in [0,T^n]; ~~~ [\alpha^n_i, \beta^n_i], i=0,1,\ldots \}.
$$ 
(When this does not cause confusion, which is in most cases, we will call it just a realization.)

\begin{definition}[Fluid sample path (FSP)]
\label{def-fsp}
Suppose, there is a fixed sequence of the process (extended) realizations on the time intervals $[0,T^n]$
\beql{eq-traj-prelim}
\{ [f^n(\cdot,t),\ell^n(t), r^n(t), \tau^n(t)], ~t\in [0,T^n]; ~~~ [\alpha^n_i, \beta^n_i], i=0,1,\ldots \},
\eeql
such that the driving sequences' realizations satisfy conditions \eqn{eq-xi-fslln}-\eqn{eq-eta-slln}.
Then, a trajectory 
\beql{eq-traj-fsp}
\{ [f(\cdot,t),\ell(t), r(t), \tau(t)], ~t\in [0,T]; ~~~ [\alpha_i, \beta_i], i=0,1,\ldots \},
\eeql
is called a {\em fluid sample path} (FSP) on the interval $[0,T]$ if the following conditions hold.

i) Points $0\le \alpha_0 \le \alpha_1 \le \ldots$ are such that there is only a finite number of them on any finite interval; $\alpha_i < \alpha_{i+1}$,
as long as $\alpha_i < \infty$; if $\alpha_i > T$, then necessarily $\alpha_i=\infty$.
Points $0\le \beta_0 \le \beta_1 \le \ldots$ are such that there is only a finite number of them on any finite interval; $\alpha_i \le \beta_i < \alpha_{i+1}$,
as long as $\beta_i < \infty$; if $\beta_i > T$, then necessarily $\beta_i=\infty$.

(ii) Functions $f(\cdot,t)$ and $\tau(t)$ are continuous in $t$. Function $[\ell(t), r(t)]$ is right-continuous with left-limits (RCLL) in $t$;
its only possible points of discontinuity are those $\alpha_i$, $i\ge 1$, that are finite. 

(iii) Trajectory \eqn{eq-traj-fsp} is a limit of trajectories \eqn{eq-traj-prelim}, as $n\to\infty$, in the following sense:

(iii.1) $T^n\to T$;

(iii.2) $\alpha^n_i\to \alpha_i$ and $\beta^n_i \to \beta_i$, for each $i$;

(iii.3) The following Skorohod space convergence holds:
$$
\{ [f^n(\cdot,t),\ell^n(t), r^n(t), \tau^n(t)], ~t\in [0,T^n]\} \rightarrow 
\{ [f(\cdot,t),\ell(t), r(t), \tau(t)], ~t\in [0,T] \}.
$$
(For the purposes of the Skorohod space metric, by convention, $[f(\cdot,t),\ell(t), r(t), \tau(t)]$ is defined for all $t\ge 0$, with its value being
constant for $t\ge T$.)
\end{definition}

Note that in Definition~\ref{def-fsp} convergence
$$
[f^n(\cdot,t),\ \tau^n(t)] \rightarrow  [f(\cdot,t), \tau(t)]
$$
is necessarily uniform on compact sets (since $[f(\cdot,t), \tau(t)]$ is continuous), and convergence 
$$
[\ell^n(t), r^n(t)] \rightarrow  [\ell(t), r(t)]
$$
is necessarily uniform on any closed bounded interval, not containing any of the points $\alpha_i$.

The following lemma describes basic properties of the FSPs, implied by the corresponding basic properties of pre-limit trajectories  and FSP definition.

\begin{lemma} 
\label{lem-fsp-prop-A-basic} 
Any FSP has the following properties.\\
(i) For any $t$, $f(x,t)$ is Lipschitz with constant $1$. Moreover, $f(r(t),t)-f(\ell(t),t)=r(t)-\ell(t)=\tau(t)$ and $f(x,t)$ is Lipschitz with constant $1/2$
in $[r(t),\ell(t)+1]$.\\
(ii) $f(x+1,t)-f(x,t) = \rho$ for any $x$ and $t$.\\
(iii) For a $t>0$, condition $\tau(t)>0$ necessarily implies that $t\in (\alpha_i, \beta_{i} \wedge T]$ for some $i$.\\
(iv) Function $\tau(t)$ is Lipschitz with constant $1$.\\
(v) Functions $\ell(t)$ and $r(t)$ are Lipschitz non-increasing, with constant $1$, on the intervals 
of the form $[\alpha_i, \beta_{i}\wedge T)$. In the sub-intervals of $[0,T)$, not intersecting with any of the intervals
$[\beta_i, \alpha_{i+1})$, $\tau(t)=0$ and  $\ell'(t)=r'(t)=0$.
\end{lemma}

{\em Proof} is rather straightforward. For example, to prove (i), note that for any $n$,$t$ and $\delta\ge 0$, $f^n(x+\delta,t)-f^n(x,t) \le \delta +2/n$;
moreover, $\delta -2/n \le f^n(x+\delta,t)-f^n(x,t) \le \delta +2/n$ if $[x,x+\delta] \in [\ell^n(t), r^n(t)]$, and
$f^n(x+\delta,t)-f^n(x,t) \le \delta/2 +2/n$ if $[x,x+\delta] \in [r^n(t),\ell^n(t)+1]$.
Given the FSP definition, this easily implies that: $f(x,t)$ is Lipschitz in $x$, with constant $1$, everywhere;
$(\partial/\partial x)f(x,t)=1$ if $x\in (\ell(t),r(t))$; $f(x,t)$ is Lipschitz in $x$, with constant $1/2$, in $[r(t),\ell(t)+1]$. 
The rest of the properties are easily established as well. We omit further details. $\Box$

For an FSP, a time point $t$ is called {\em regular}, if the derivatives $\tau'(t), \ell'(t), r'(t)$ exist. In view of Lemma~\ref{lem-fsp-prop-A-basic}(iv)-(v), {\em almost all points (w.r.t. Lebesgue measure) are regular.}

We see that an FSP describes the evolution of the distribution of the continuous ``particle mass,'' or ``fluid,'' obtained as a limit of rescaled pre-limit particle distributions. For an FSP we will also use the natural notion of a cluster: we say that a cluster $[\ell(t),r(t)]$,
of length $\tau(t)=r(t)-\ell(t)$,
 exists at time $t$,
if $t\in [\alpha_i, \beta_{i})$ for some $i$; in this case the rest of the unit circle, that is interval $[r(t),\ell(t)+1]$ is a sparse interval. 
A sub-interval of a sparse interval we will call a sparse sub-interval.
Note that it is possible that a cluster in an FSP exists even when it has zero length, $\tau(t)=0$. Of course, if
$\tau(t)>0$ then a cluster necessarily exists. If there is no cluster at time $t$, then necessarily $\tau(t)=0$ and the entire unit circle $[0,1]$
is a sparse interval. The total particle mass $f(x+1,t)-f(x,t) = \rho$ is constant at all times.
The particle mass density $(\partial/\partial x) f(x,t)$ within a cluster is exactly $1$, and the density within a sparse interval is at most $1/2$.

\begin{figure}[h!]
\centering
\includegraphics[width=1.0\linewidth]{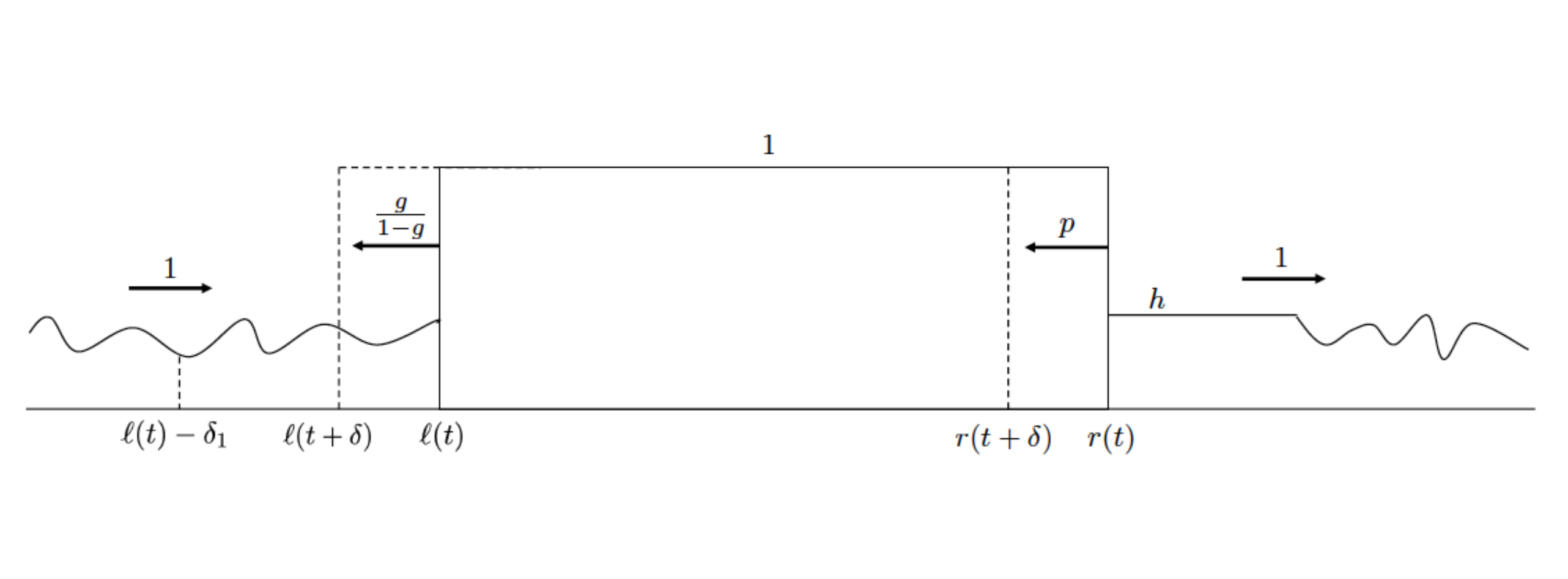}
\caption{Cluster dynamics illustration}
\label{fig:cluster-dyn}
\end{figure}

The following Lemma~\ref{lem-fsp-prop-A-cluster-deriv} describes the basic FSP dynamics. We precede it with an informal description of this dynamics. (See Figure~\ref{fig:cluster-dyn} for an illustration.) The particle mass outside a cluster (if any) simply moves to the right at the constant speed $1$ until and unless it ``hits'' the cluster.
In a time interval when a cluster $[\ell(t),r(t)]$ exists,  its right end $r(t)$ moves left at the constant speed $p$. (Because this is the rate at which particles in a pre-limit system ``break-away'' from the cluster on the right.) Moreover, as the right end of a cluster moves left, it leaves a sparse sub-interval of the density exactly $h=p/(1+p)$ in its ``wake." (Because, in pre-limit process, the average distance between two consecutive breaking away particles is $1+1/p = 1/h$.) Specifically, in a time interval $[t,t+\delta]$, a sparse sub-interval $[r(t)-\delta p, r(t) + \delta]$ of density $h$ is created, where $r(t+\delta)= r(t)-\delta p$ is the right-end location at time $t+\delta$, and $r(t)+\delta$ is how far the right edge of the sparse sub-interval moved right by time $t+\delta$. Now let us consider how the left end $\ell(t)$ moves. Consider a small interval $[t,t+\delta]$. 
Let $[\ell(t)-\delta_1,\ell(t)]$ be the sparse sub-interval immediately to the left of the cluster at time $t$, which hits -- and merges into -- the cluster in time interval $[t,t+\delta]$.
Then the $\delta_1$ above is the unique solution of equation
\beql{eq-time-delta1}
 \delta=\delta_1 - [f(\ell(t),t)-f(\ell(t)-\delta_1,t)].
\eeql
Indeed, the time $\delta$ for the sparse sub-interval $[\ell(t)-\delta_1,\ell(t)]$ (at $t$) to join the cluster is exactly equal to the ``empty space'' 
in $[\ell(t)-\delta_1,\ell(t)]$ (at $t$), which is $\delta_1 - [f(\ell(t),t)-f(\ell(t)-\delta_1,t)]$; solution $\delta_1$ is unique, because the RHS of
\eqn{eq-time-delta1}, as function of $\delta_1$, is Lipschitz increasing with derivative $\ge 1/2$ (recall that the density within a sparse interval is at most $1/2$).

The following lemma describes the above basic dynamics of an FSP formally. Its proof is also fairly straightforward -- 
most of it uses the corresponding properties of pre-limit trajectories  and FSP definition.
We will only provide the key points and comments, omitting details. 

\begin{lemma} 
\label{lem-fsp-prop-A-cluster-deriv}
Any FSP has the following properties.

(i)  Consider a  time point $t$ such that $t\in (\alpha_i, \beta_{i} \wedge T)$.
 (Recall that for a $0<t<T$, $\tau(t) > 0$  necessarily implies $t\in (\alpha_i, \beta_{i} \wedge T)$.)
 For any sufficiently small time $\delta>0$, the following holds. Denote by $\delta_1$ the unique positive solution to
 equation 
 \beql{eq-time-delta}
 \delta=\delta_1 - [f(\ell(t),t)-f(\ell(t)-\delta_1,t)].
 \eeql
   Then, the FSP state at time $t+\delta$ is as follows:
 \beql{eq-ell-moves}
 \ell(t+\delta) = \ell(t) - (\delta_1-\delta);
 \eeql
 \beql{eq-r-moves}
 r(t+\delta) = r(t) - p \delta;
 \eeql
  for $x \not\in (\ell(t+\delta) , r(t + \delta))$ (this is the part of sparse interval that just shifted, without any interaction with the cluster) we have
 $$
 f(x,t+\delta)= f(x-\delta,t);
 $$
 for $x\in (\ell(t+\delta), r(t+\delta))$ (within cluster at time $t+\delta$)
 $$
 (\partial/\partial x) f(x,t+\delta)= 1;
 $$
 for $x \in (r(t+\delta), r(t) + \delta)$  (the ``wake'' of right end of cluster)
 $$
(\partial/\partial x) f(x,t+\delta)= h.
 $$
 In addition, if this $t$ is regular,
 \beql{eq-l-deriv}
 \ell'(t) = - \frac{g}{1-g}, ~~\mbox{where}~g=g(t)= \frac{\partial^-}{\partial x} f(\ell(t),t)
 \eeql
   and $r'(t)=-p$. 
   
 (ii)  If $t$ is regular, then $\tau'(t)>0$  if and only if $g>h$.
 
 (iii) If $t$ is regular, $t\in (\alpha_i, \beta_{i}\wedge T)$ and $\tau(t)=0$, then, necessarily,
 $\tau'(t) =0$ and $g=h$.
 
 (iv) For any  interval $[t,t+\delta] \subset [0,T)$ when a cluster does not exist (i.e. not intersecting with any of the intervals $[\alpha_i, \beta_{i})$),
the particle density moves to the right at speed $1$:
 $$
  f(x,t+\delta) = f(x-\delta,t), ~\ell(t+\delta)=r(t+\delta)=\ell(t)+\delta=r(t)+\delta, ~\tau(t)=0.
 $$
 
(v) Denote by $\mu(t)$ the total Lebesgue measure of all points $x\in [0,1)$ 
 (within the circle) at time $t$,  which are outside the cluster $[\ell(t),r(t)]$ (if any) and where the density 
 $$
 \frac{\partial}{\partial x} f(x,t) > h.
 $$
Function $\mu(t)$ is Lipschitz (and then the derivative exists for almost all $t$), and given the properties (i)-(iv), it is non-increasing.
Then, for almost all $t>0$ (w.r.t. Lebesgue measure),
either
 $$
\mu'(t) = 0,
 $$
 or this time $t$ is regular, $\tau'(t)>0$ (and then $\tau(t)>0$, $t\in (\alpha_i, \beta_{i}\wedge T)$), and
 $$
g> h ~~\mbox{and}~~ \mu'(t) = - 1 + \ell'(t) = - \frac{1}{1-g} \le - \frac{1}{1-h} = -  (1+p).
 $$

 (vi) The Lebesgue measure of the time points $t$ where $\tau'(t)>0$ is upper bounded
 by $1/(1+p)$.
 \end{lemma}
 
 From now on we make the definition of a regular point 
 more restrictive by requiring, in addition, that the conclusions of Lemma~\ref{lem-fsp-prop-A-cluster-deriv}(v) hold.
 Still, almost all time points $t$ are regular.
  
 {\em Proof.} (i) The meaning of \eqn{eq-time-delta} (which is same as \eqn{eq-time-delta1}) is explained in the informal definition above.
 Consider a fixed $\delta>0$ and time $t$. For each pre-limit trajectory, as in  FSP definition, 
 define $\delta_1^{(n)}\ge 0$ as the smallest number such that all particles located at time $t$ in the interval $[\ell^n(t)-\delta_1^{(n)},\ell^n(t)]$,
 will join the cluster by time $t+\delta$. It is easy to observe that
 $$
 \delta=\delta_1^{(n)} - [f(\ell^n(t),t)-f(\ell^n(t)-\delta_1^{(n)},t)] + O(1/n).
$$
It remains to consider the limit in $n$ (as in the FSP definition) and observe that $\delta_1^{(n)}$ must converge exactly to the $\delta$ as defined.
All properties before \eqn{eq-l-deriv} then easily follow. Property \eqn{eq-l-deriv} is a differential form of \eqn{eq-time-delta}. Namely,
we use the fact that $t$ is regular and therefore $\ell'(t)$ exists; let $\delta\downarrow 0$; from \eqn{eq-ell-moves} we have
$$
\lim_{\delta\to 0} \delta_1/\delta = 1- \ell'(t);
$$
dividing \eqn{eq-time-delta} by $\delta$, taking $\delta\downarrow 0$ limit and substituting the above display, we obtain \eqn{eq-l-deriv}.
The last property is from \eqn{eq-r-moves}.

(ii) and (iii) follow from (i). (iv) is obvious.

(v) Obviously, $\mu'(t)=0$ for almost all time points $t$ such that a cluster does not exist. Consider now a $t$ which is strictly within a time interval when a cluster exists. Once again, consider the movement of the left end $\ell$ of the cluster in a small interval $[t,t+\delta]$. Then,
using \eqn{eq-l-deriv},
\eqn{eq-time-delta} can be written as 
\beql{eq-measure}
- \delta_1 = \int_t^{t+\delta} [-1+\ell'(s)]ds = \int_t^{t+\delta} \frac{-1}{1-g(s)}ds,
\eeql
where $-1+\ell'(s)$ is the instantaneous (negative) rate at which the sparse sub-interval $[\ell(t)-\delta, \ell(t)]$ of length $\delta_1$, which existed
at time $t$, ``shrinks'' due to its right end $\ell(s)$ moving at speed $\ell'(s)$ left and its left end $\ell(t)-\delta_1+(s-t)$ moving at speed $1$ right.
The integrand in the RHS exists for almost all $t$ (and is unique up the time subsets of zero Lebesgue measure).
Note that, for a given $t$, the density $(\partial/\partial x)f(x,t)$ exists a.e. in $x$, and therefore $(\partial^-/\partial x)f(x,t)=(\partial/\partial x)f(x,t)$ a.e. (Recall that $g(t)$ is defined as left derivative.)
Analogously to \eqn{eq-measure}, it is easy to obtain its generalization, which gives the (negative) increment of $\mu$ in the interval
$[t,t+\delta]$, by integrating over only those times $s$, where $(\partial^-/\partial x) f(\ell(s),s) = g(s) > h$: 
\beql{eq-measure-mu}
\mu(t+\delta) - \mu(t) = \int_t^{t+\delta} \frac{-1}{1-g(s)} \I\{g(s) > h\} ds.
\eeql
Therefore, for almost all $s\in [t,t+\delta]$,
$$
\mu'(s) =  \frac{-1}{1-g(s)} \I\{g(s) > h\}.
$$
From here, along the fact that almost all times $s$ are regular, all properties stated in (v) follow.

(vi) Follows from (v).
 $\Box$

  \begin{lemma} 
\label{lem-fsp-prop-A-low}
Suppose, $\rho < h$. Then any FSP is such that the following properties hold.

(i) For any $i$ such that $\alpha_i < \infty$, we have $\beta_i \wedge T - \alpha_i < 1$. \\
Consequently, for any $t_0 \le T -1$ there exist $t_0 < t_1 \le t_0 +1$ such that $\tau(t_1)=0$.

(ii) 
\beql{eq-tau-int-bound}
\int_0^T \tau(t) dt \le \bar C = 1 + \frac{1}{2(1+p)}.
\eeql
\end{lemma}

{\em Proof.}  (i) Suppose not, i.e. $(t_0,t_0+1]$ is entirely within an interval of the form $[\alpha_i, \beta_i \wedge T)$.
Then,
at time $t_0+1$ the FSP state is such that everywhere outside the cluster the density 
is equal to $h$. (This follows from Lemma~\ref{lem-fsp-prop-A-cluster-deriv}(i), density $h$.)  
This is impossible, because it would imply that $\rho \ge h$.

(ii) Consider the set of points $t\in [0,T)$, where $\tau(t)>0$. This set consists of possibly an interval $[0, s)$, whose length $s$ cannot exceed $1$ (by (i)), and a countable number of open intervals $(t_1,t_2)$ with lengths also not exceeding $1$ (by (i)). Obviously,
$$
\int_0^s \tau(t) dt \le s \le 1.
$$ 
Consider any of the open intervals $(t_1,t_2)$. We have 
$$
\max_{t\in (t_1,t_2)} \tau(t) \le \int_{t\in (t_1,t_2):~\tau'(t) > 0} \tau'(t) dt \le \int_{t\in (t_1,t_2):~\tau'(t) > 0} [-\ell'(t)] dt
\le \frac{1}{2} \mathcal L \{t\in (t_1,t_2):~\tau'(t) > 0\},
$$
and therefore (recall that $t_2-t_1\le 1$) 
$$
\int_{t_1}^{t_2} \tau(t) dt \le \frac{1}{2} \mathcal L \{t\in (t_1,t_2):~\tau'(t) > 0\}.
$$
Summing up over all intervals where $\tau(t)>0$, we obtain
$$
\int_0^T \tau(t) dt \le 1 + \frac{1}{2} \mathcal L \{t\in (0,T):~\tau'(t) > 0\} \le 1 + \frac{1}{2} \frac{1}{1+p},
$$
where in the last inequality we used Lemma~\ref{lem-fsp-prop-A-cluster-deriv}(vi).
$\Box$

Suppose $h < \rho < 1$.  Recall $h=p/(1+p)$. State $[f^*(\cdot),\ell^*, r^*, \tau^*]$ is called a {\em main equilibrium state} (MES),
if it satisfies the following conditions: $\tau^* = \frac{\rho-h}{1-h}$, $r^* = \ell^* + \tau^*$, $f'(x)=1$ for $x\in (\ell^*,r^*)$,
 $f'(x)=h$ for $x\in (r^*, \ell^*+1)$. The definition shows that a MES is essentially unique,  up to a space shift of the cluster; 
 therefore, we will often refer to any MES as {\em the} MES. Note that if the steady-state of a system with large $n$ is close to MES, then
 the flux is close to $\phi^*=\phi^*(\rho) \doteq (1-\tau^*)h = p(1-\rho)$.
 
  \begin{lemma} 
\label{lem-fsp-equil-states}
Suppose, $\rho > h$. Then any FSP is such that the following properties hold.

 (i) If FSP initial state is a MES, the FSP is unique and is stationary, i.e. it stays in this state (up to shift of the cluster).
 
 (ii) If $\alpha_0 =0$ and $\beta_0 \ge 1-\rho$, 
  then the FSP state at time $1-\rho$ is a MES, and the trajectory is stationary (staying in a MES) in the interval $[1-\rho,T]$.
\end{lemma}

{\em Proof.} (i) According to derivatives' expressions for a MES, it cannot change with time, up to the cluster moving left at the constant speed $p$.

(ii) The entire particle mass, which is originally outside the cluster is at most $\rho$.
We have $\beta_0 \ge 1-\rho$, so the cluster exists for the time at least $1-\rho$. The time interval $[0,1-\rho]$ is long enough for all the mass which was originally outside the cluster to join the cluster. Therefore, at some time point in $t_0 \in [0,1-\rho]$, a state is reached, consisting of a cluster and a sparse interval with constant density $h$ (by Lemma~\ref{lem-fsp-prop-A-cluster-deriv}(i).)
This is equivalent to the state being a MES. The same argument as in the proof of (i) shows that the trajectory must stay in MES
in  $[t_0,T]$.
$\Box$

\begin{lemma} 
\label{lem-fsp-conv-A}
Consider a sequence of systems, with $n\to\infty$, with arbitrary random initial states $[f^n(\cdot,0),\ell^n(0), r^n(0), \tau^n(0)]$.
For each $n$ consider some stopping time $T^n$ (finite or infinite). 
Then, with prob. 1, any subsequence of (extended) realizations has a further subsequence, along which 
the convergence to an FSP holds in the sense of Definition~\ref{def-fsp}(iii).
\end{lemma}

{\em Proof.} This type of a fluid-limit result is standard. If the sequence of processes is constructed on the common probability space, as specified above, the driving sequences satisfy the properties required in the FSP definition w.p.1. 
Then, w.p.1., from any subsequence we can choose a further subsequence,
such that there is a convergence to some trajectory. The latter trajectory must be an FSP, by the FSP definition. 
We omit further details.
$\Box$

\begin{lemma} 
\label{lem-fsp-conv-A-high-stopping}
Suppose $\rho > h$. Consider a sequence of systems, with $n\to\infty$, with arbitrary random initial states $[f^n(\cdot,0),\ell^n(0), r^n(0), \tau^n(0)]$. For a fixed constant $B\ge 1-\rho$ and each $n$ consider a stopping time $T^n$, which is the minimum of constant $B$ and the first time after $0$, when the process is in an ideal state. 
Then, w.p.1, any subsequence of (extended) realizations has a further subsequence, along which it converges to 
an FSP of duration $T\le B$. This FSP satisfies one of the following two properties:\\
(a) $T < 1-\rho$ and $\tau(T)=0$;\\
(b) $T=B$ and the FSP state in $[1-\rho, B]$ is a MES.
\end{lemma}

{\em Proof.} The probability 1 subsequential convergence to FSPs is due to Lemma~\ref{lem-fsp-conv-A}. Then the required FSP properties follow from Lemma~\ref{lem-fsp-equil-states}(ii).
$\Box$

\subsection{A version of the ballot theorem}
\label{sec-ballot-prop}

The following combinatorial fact plays a key role in the proofs of our main results for the high density, $h < \rho < 1/2$. We will derive it 
from Theorem 1.2.5 -- a ballot theorem -- in \cite{Tac1967}.

\begin{proposition}
\label{th-ballot}
Suppose an integer $n\ge 1$ and a sequence of real non-negative numbers, $k_1, k_2, \ldots, k_n$ are fixed. Let us extend the definition of the sequence $\{k_r\}$ to all integer $r\ge 1$ by periodicity, $k_{r+n}=k_r$, and denote
$\psi(j)=  \sum_{\ell=1}^{j} k_{\ell}$ for $j\ge 0$, where $\psi(0)=0$ by convention. Suppose a real number $m >\psi(n)$ is fixed.
Let $N$ be the number of those indices $j \in \{0,1,\ldots,n-1\}$, for which
\beql{eq-bal1}
\psi(j+r) -\psi(j) < \frac{m}{n} r, ~~~  r=1,\ldots,n.
\eeql
Then,
\beql{eq-bal2}
N \ge \lceil n\left( 1- \frac{\psi(n)}{m}\right)   \rceil.
\eeql
\end{proposition}

The proof is in Appendix ~\ref{app-ballot}.

\subsection{A-perturbations model: Main results}
\label{sec-A-pert-main-res}

\begin{theorem}
\label{th-A-low} 
Consider the model with A-perturbations. Assume $\rho < h$. 
Then, as $n\to\infty$, the sequence of the stationary distributions of
 $[f^n(\cdot,t),\ell^n(t), r^n(t), \tau^n(t)]$ is such that its any subsequential weak limit is concentrated on the 
 states with $\tau=0$. 
In other words, as $n\to\infty$, $\tau^n(\infty) \Rightarrow 0$ [equivalently, $\E \tau^n(\infty) \to 0$] 
and, consequently, 
the limiting flux (typical flux) $\phi(\rho)= \rho$.
 \end{theorem}

{\em Proof of Theorem~\ref{th-A-low}.}
It suffices to show that any distributional limit is such that $\E \tau =0$. Suppose not. Consider a fixed interval $[0,B]$, and stationary versions of the processes in this interval. Then there exists $\epsilon>0$ such that, for any $B>0$,
\beql{eq-contr1}
\liminf_{n\to\infty} (1/B) \E \int_0^B \tau^n(t) dt \ge \epsilon > 0.
\eeql  
However, using our construction above,
we can construct all these (stationary) processes on a common probability space, so that, w.p.1.,
we have a subsequential convergence to an FSP. Convergence to an FSP, in particular, means that $\tau^n(\cdot)\to \tau(\cdot)$ u.o.c..
Then, for any $B>0$, from Lemma~\ref{lem-fsp-prop-A-low}(ii) we obtain 
$$
\limsup_{n\to\infty} \E \int_0^B \tau^n(t) dt \le \bar C.
$$
This contradicts the fact that \eqn{eq-contr1} must hold for all $B>0$.
$\Box$

\begin{lemma}
\label{lem-ballot-appl}
Consider the model with A-perturbations. Assume $h < \rho < 1/2$. Consider the system without time/space rescaling.
Consider a (random) initial system state, at initial time $0$, which
is formed as follows: we pick an arbitrary ideal state and apply a perturbation to it.
(In other words, a particle is picked uniformly at random, and it is placed into one of the empty sites, chosen uniformly at random.)
Then, there exists $\delta=\delta(p,\rho)>0$ such that, uniformly in all the original (pre-perturbation) ideal states and uniformly in all sufficiently large $n$ (and then in all $n$),
the following event occurs with probability at least $\delta$: The perturbation creates a cluster and this cluster will {\em not} disappear within (unscaled) time interval $[0,(1-\rho)n-4]$.
\end{lemma}

{\em Proof.} It will be convenient to relabel time, so that the initial time is $1$ (instead of $0$).
Consider the initial state after the perturbation. Consider the relocated particle, which we will refer to as a ``seed'' particle; the site to which it relocates we will call the ``seed'' site. If the perturbation happens to form a cluster, let us consider the dynamics of this cluster, starting the initial time $1$. (If the perturbation did not form a cluster, let us view such event as a cluster ``disappearance'' immediately at time $1$.) Obviously, as long as the cluster exists, the particles with leave it on the right as a Bernoulli process, with ``success'' (leaving) probability $p$, starting time $2$. (By convention, a particle does {\em not} leave the cluster at initial time $1$.) Of course, the process of ``failures'' (non-departures) is Bernoulli with probability $1-p$. Fix $\epsilon=(\hat p -p)/3$, where $\hat p = \rho/(1-\rho) > p$, because $p=h/(1-h)$.
(The meaning of $\hat p$ will be explained shortly.) Let $D(t), ~t=1,2,\ldots,$ be the cumulative number of ``successes'' of the Bernoulli process described above (with $D(1)=0$ by convention), by and including time $t$, and let $\bar D(t) = t-D(t)$ be the corresponding cumulative number of ``failures'' 
(with $\bar D(1)=1$ by convention).
We observe that 
\beql{eq-bound-right}
\delta_1 = \P\{\bar D(t) \ge [(1-p) - \epsilon]t, ~\forall t\ge 1\} > 0.
\eeql

Now consider how the cluster (if any) ``grows'' on the left due to particles joining it on the left. (By convention,
if the seed particle finds another one immediately to the left of the seed site, 
we assume that the latter particle joined the cluster immediately at the initial time $1$.)
Note that the initial cluster (if any) has at most 1 particle to the right of the seed particle and, therefore, at time $1$, at least $\rho n -3$ particles (i.e., ``almost all'') form a sparse interval to the left of the cluster. Therefore, at least $\rho n -2$ particle which are initially to the left of the seed particle will move unobstructed, at speed $1$, until and unless they join the cluster.

Let us define the following process over the time interval $1,2, \ldots, n - (\rho n -1) = (1-\rho)n +1$. (It is such that it exactly describes the process $A(t)$ of the cumulative number of particles joining the cluster from the left, in the (a bit shorter) time interval $1,2,\ldots, n -(\rho n -1) - 2\cdot 2 = (1-\rho)n -3$.) We will illustrate the definition by an example. Let $n=10$ and the total number of particles is $\rho n=4$. Consider the configuration of the $\rho n-1=3$ particles on the circle, {\em excluding the seed particle}, with respect to the seed site. Starting from the site immediately to the left from the seed site, and moving left, the configuration of the particles (excluding seed particle) is, for example, this sequence of length $n=10$:
\beql{eq-seq1}
1,0,0,1,0,0,0,0,1,0
\eeql
(Note that the last element of this sequence is necessarily $0$, because it corresponds to the seed site, and we do {\em not} include the seed particle into this sequence.) Consider now the sequence of length $n - (\rho n -1) = 10 - 3 = 7$, obtained from \eqn{eq-seq1} by removing  each $0$, which immediately follows a $1$. We obtain
\beql{eq-seq2}
1,0,1,0,0,0,1
\eeql
The cumulative number of $1$'s by time (= place in the sequence) $t$ is the function $A(t), ~t=1,2, \ldots, (1-\rho)n +1$. 
It is easy to see that it {\em would} exactly describe the number of particles joining the cluster on the left, assuming the cluster would continue to exist, and {\em assuming all particles except the seed one would move freely at rate $1$}. It is also easy to observe that, in a little shorter interval
$t=1,2, \ldots, (1-\rho)n -3$, $A(t)$ is {\em exactly} the number of particles joining the cluster on the left, assuming the cluster would continue to exist. 

Clearly, the process $A(t)$ satisfies the following conservation law: $A((1-\rho)n +1) = \rho n -1$.  
Let us denote by $\bar A(t) = t- A(t)$ the corresponding cumulative process of ``failures'' (zeros). The conservation law for $\bar A(t)$
is: $\bar A((1-\rho)n +1) = ((1-\rho)n +1) - (\rho n -1) = (1- 2\rho)n + 2$. Note that the average slope of $\bar A(t)$ is
\beql{eq-Aslope}
\frac{(1- 2\rho)n + 2}{(1-\rho)n +1} = \frac{(1- 2\rho) + 2/n}{(1-\rho) +1/n},
\eeql
and is close to 
$$
\frac{1- 2\rho}{1-\rho} = 1-\hat p, ~~  \hat p = \rho/(1-\rho).
$$
Therefore, for all large $n$, the average slope \eqn{eq-Aslope} is at most $(1-\hat p) + \epsilon < (1-\hat p) + 2\epsilon = (1- p) - \epsilon$.
Notice that, given locations of all particles except the seed particle on the circle, since the seed particle chooses one of the empty sites (to become the seed site) uniformly at random,
all cyclical permutations of a sequence \eqn{eq-seq2} appear with equal probabilities. Using the latter observation, we can apply Proposition~\ref{th-ballot}, to obtain the following estimate: for all sufficiently large $n$,
\beql{eq-bound-left}
\P\{\bar A(t) < [(1-\hat p) + 2\epsilon]t, ~\forall~1 \le t \le (1-\rho)n +1\} \ge \delta_2 = \frac{\epsilon}{(1-\hat p) + 2\epsilon} > 0.
\eeql

Combining estimates \eqn{eq-bound-right} and \eqn{eq-bound-left}, we conclude that, uniformly on all initial (pre-perturbation) ideal states and in all sufficiently large $n$, with probability at least $\delta = \delta_1 \delta_2$, a perturbation will create a cluster at time $1$ and this cluster will not disappear by and including time $(1-\rho)n -3$.
$\Box$

As a corollary of Lemmas~\ref{lem-fsp-conv-A-high-stopping} and \ref{lem-ballot-appl}, we obtain the following

\begin{lemma}
\label{lem-ballot-appl2}
Consider the model with A-perturbations. Assume $h < \rho < 1/2$. For each $n$ consider the scaled process, with (random) initial 
state formed  by a perturbation of an arbitrary ideal state. Fix $B\ge 1-\rho$. Consider the stopping time $T^n$ which is the minimum of $B$ and the first time the process hits an ideal state. Then any subsequential (in $n\to\infty$) weak limit of the state distributions at time $T^n$
is such that, with probability at least $\delta$ (defined in Lemma~\ref{lem-ballot-appl}) the state is a MES.
(The subsequential weak limits exist, because the entire state space of $[f^n(\cdot,t),\ell^n(t), r^n(t), \tau^n(t)]$ is compact, if we consider the unique ``standard'' version of each state, satisfying $f^n(0,t)=0$ and $\ell^n(t)\in [0,1)$.)
\end{lemma}

\begin{theorem}
\label{th-A-high}
Consider the model with A-perturbations. Assume $h < \rho < 1/2$.  Then, as $n\to\infty$, any subsequential weak limit of the sequence of the stationary distributions of
 $[f^n(\cdot,t),\ell^n(t), r^n(t), \tau^n(t)]$ is a distribution concentrated on the main equilibrium states.
 Consequently, the limiting flux (typical flux) $\phi(\rho)= \phi^*= (1-\tau^*)h = p(1-\rho)$.
  \end{theorem}
  
  \begin{remark} Theorem~\ref{th-A-high} actually holds for $h < \rho < 1$. The proof for 
 $1/2 < \rho \le 1$ is a simplified version of that for $h < \rho < 1/2$, because when $\rho > 1/2$ a cluster always exists and A-perturbations never occur. If $\rho = 1/2$, the corresponding sequence of pre-limit systems is such that an ideal state either exists infinitely often, or does not exist infinitely often, or both. In any case, if we consider the corresponding subsequences of systems, 
 the proof is the same as for either $h< \rho < 1/2$ or $h> 1/2$.
  \end{remark}

{\em Proof of Theorem~\ref{th-A-high}.} 
For each $n$, consider the process with any fixed initial state at $0$. Consider 
stopping time $\theta^n_1$, which is the minimum of 2 and the first time (after 0) the system is in an ideal state. 
Consider the random system state at $\theta^n_1$. Consider its any subsequential limit in distribution. 
By Lemma~\ref{lem-fsp-conv-A-high-stopping}, this distributional limit is such that, almost surely, the state is either a main equilibrium state (MES) or it is an ideal state. (A small subtlety here: we use the fact that, as $n\to\infty$, the probability that a perturbation occur immediately at a time of hitting an ideal state, vanishes.) 

Consider stopping time $\theta^n_2$, which is the minimum of $\theta^n_1+2$ and the first time after $\theta^n_1$ the system is in an ideal state. Note that, as $n\to\infty$, the time to wait for a perturbation converges to exponentially distributed one with mean $1/\lambda$;
in particular, uniformly in all large $n$, the probability that a time to wait for a perturbation is less than $1$, is at least 
$\epsilon = (1-e^{-\lambda})/2$. 
If the state at $\theta^n_1$ was ideal, then with probability $\ge \epsilon$, there will be a perturbation in $[\theta^n_1, \theta^n_1+1]$.
Using Lemma~\ref{lem-ballot-appl2}, we conclude that any subsequential limit in distribution of the system state at $\theta^n_2$ is
either a MES {\em with probability at least $\zeta=\epsilon \delta>0$}, or it is an ideal state. (The $\delta$ is as defined in 
Lemmas~\ref{lem-ballot-appl} and \ref{lem-ballot-appl2}.)

We define stopping times $\theta^n_k$ for all $k\ge 2$ analogously to the definition of $\theta^n_2$. By induction,
we conclude that any distributional limit of the system state at $\theta^n_k$
is such that it is either a MES {\em with probability at least $1-(1-\zeta)^{k-1}$}, or it is an ideal state.
Note that $\theta^n_k \le 2k$. 
Finally, using Lemma~\ref{lem-fsp-equil-states}(i), we establish
that any subsequential limit in distribution of the system state at time $2k$ is such that 
with probability at least $1-(1-\zeta)^{k-1}$, the state is a MES. 
This is true for {\em any} fixed initial states for each $n$, and any $k$.
This implies that the limit of stationary distributions of the process is concentrated on MESes.
$\Box$

\section{I-perturbations model}
\label{sec-i-perturb}

Unless specified otherwise, in this section we consider the model with I-perturbations.

\subsection{Some preliminary facts.}
\label{sec-i-prelim}

As a corollary of Lemma~\ref{lem-simple}, we obtain the following fact.

\begin{lemma}
\label{lem-simple-low}
Suppose $0< \rho <  h$. Assume there are {\em no perturbations} in the (scaled) time interval $[0,1]$. Then, uniformly on the initial states at time $0$, which may contain one or multiple clusters,
$$
\lim_{n\to\infty} \P\{\mbox{At (scaled) time $1$ the system is in an ideal state}\} = 1.
$$
\end{lemma}

{\em Proof of  Lemma~\ref{lem-simple-low}.} For each $n$, consider an arbitrary system state at time $0$, and pick any cluster (if any). Let us show that, uniformly on the initial states and the choice of the cluster,
\beql{eq-dissap}
\lim_{n\to\infty} \P\{\mbox{The cluster will disappear before (scaled) time $1$}\} = 1.
\eeql
WLOG, let the right edge of the cluster be located at (scaled) coordinate 0. Then, the ``best case,'' in terms of the cluster {\em not} disappearing for as long as possible, is when this cluster contains all particles in the system. 
From the law of large numbers, we easily verify that, for a fixed $\epsilon >0$ such that $\rho+\epsilon < 1$,
as $n\to\infty$, 
$$
\lim_{n\to\infty} \P\{\mbox{The cluster ``dissolves'' before the particle initially at $0$ reaches point $-\rho-\epsilon$}\} = 1.
$$
This implies \eqn{eq-dissap}.

We are not done, because the total number of initial clusters may increase to infinity with $n$. 
However, the following observation resolves this potential issue. Fix a small $\epsilon >0$, such that $\rho+\epsilon<1$.
The argument we used above for a single cluster easily generalizes to show the following property: 
\beql{eq-dissap-gen}
\lim_{n\to\infty} \P\{\mbox{All clusters, which at time $0$ overlap with $[-\epsilon,0]$, will disappear before (scaled) time $1$}\} = 1.
\eeql
Indeed, if at time $0$ we ``fill in'' all holes within the scaled interval $[-\epsilon,0]$ with additional particles, we produce a single cluster (completely covering $[-\epsilon,0]$), then the original and modified processes can be coupled so that the disappearance of the cluster 
in the modified system implies the disappearance of {\em all} clusters initially overlapping with $[-\epsilon,0]$ in the original system. Note that
the total (scaled) number of particles in the modified system is at most $\rho+\epsilon$. Therefore, \eqn{eq-dissap} applies to the modified system, which proves \eqn{eq-dissap-gen}. Obviously, \eqn{eq-dissap-gen} implies the lemma statement, because the circle can be covered by 
a finite number of $\epsilon$-long intervals.
$\Box$

From Lemma~\ref{lem-simple-low} we obtain the following corollary. 

\begin{lemma}
\label{lem-partial-low}
Consider the system with I-perturbations. 
Suppose $\rho <  h$. Then, uniformly on all
sufficiently large $n$ and all initial process states, the probability that the process reaches an ideal state at some time $t_1 < 1$ is at least
$e^{-\lambda}/2$.
\end{lemma}

{\em Proof.}
The probability that there are no perturbations in $[0,1]$ converges to $e^{-\lambda}$. It remains to apply Lemma~\ref{lem-simple-low}.
$\Box$

\begin{lemma}
\label{lem-finite-num-clusters}
Suppose $\rho <  h$. Then, uniformly in $n$, the steady-state number of clusters is stochastically upper bounded by a proper (finite w.p.1) random variable.
\end{lemma}

{\em Proof.} It suffices to show that the lemma statement holds for all sufficiently large $n$.
Fix an arbitrary integer $k$. 
By Lemma~\ref{lem-partial-low}, uniformly on all initial states at time $0$, the probability that the process reaches an ideal state in $[0,k]$
is at least $1- (1-e^{-\lambda}/2)^k$. Considering the time interval from hitting an ideal state within $[0,k]$ until $k$, and taking into account the (Bernoulli) structure of the perturbation process, we easily obtain the following. (Recall also that a perturbation can increase the number of clusters by at most $2$.) Uniformly on the initial states and in all large $n$, with probability at least $1- (1-e^{-\lambda}/2)^k$, the number
of clusters in the system is stochastically dominated by a random variable $2(1+H_{2\lambda k})$, where $H_{2\lambda k}$ has Poisson distribution with mean $2\lambda k$. Since this is true for all integers $k$, we obtain a uniform (in $n$) stochastic upper bound on the 
number of clusters in steady-state.
$\Box$

\subsection{FSPs and their properties}
\label{sec-i-perturb-fl}

In this subsection we define FSPs for I-perturbation process, describe their properties, and give corresponding fluid limit results.
For our purposes, it will suffice to consider FSPs with a {\em finite number of initial clusters}. 
The definitions and results are described informally, because, in essence, they are straightforward generalizations 
of those for the A-perturbation model. We believe this informal description is sufficient -- an interested reader can easily fill in all formalities.

FSP on a (finite or infinite) interval $[0,T]$ has the following structure:
\beql{eq-traj-fsp-i}
\{ f(\cdot,t), ~t\in [0,T]; ~~~~[\ell_i(t), r_i(t), \tau_i(t)], t\in [0,T], ~ [\alpha_i, \beta_i], ~~ i=0,1,\ldots \}.
\eeql
Function $f(\cdot,t)$ describes the FSP state at $t$. The times $\alpha_i \le \beta_i \le \infty$ are the beginning and end times of the existence of $i$-cluster, and
$\ell_i(t), r_i(t), \tau_i(t)=r_i(t)-\ell_i(t)$ are the left-end, right-end and the length of the $i$-th cluster at time $t$. There is at most a finite number of initial clusters (at time $0$) -- they are indexed by $i=0,1, \ldots$ in an arbitrary order; for those initial clusters,
$\alpha_i =0$. FSP has at most a finite number of clusters by any finite time $t$ (i.e., those with $\alpha_i \le t$). At any given time $t$, the clusters present at $t$ do not overlap. New clusters are indexed in the order of their appearance, i.e. in the order of increasing $\alpha_i$. 
We adopt the following convention: if in a pre-limit process a cluster divides into two due to a perturbation, we ``ignore'' this fact and still treat the divided cluster as one; this does not cause problems, because the ``holes'' created within a cluster remain ``invisible'' in the corresponding FSP {\em on any finite time interval}, in the sense that their ``size'' remains $0$ and their existence does not affect the dynamics (derivatives) of the left and right ends, $\ell_i(t)$ and $r_i(t)$, of the original cluster in the FSP. Given this convention, except possibly time $0$, no more than one new cluster may appear 
(due to  a perturbation) at any time $t>0$. 
By another convention, before cluster $i$ exists, i.e. for $t< \alpha_i$, $\ell_i(t)=r_i(t)=\tau_i(t)=0$; and after it disappears, i.e. for $t\ge \beta_i$, $\tau_i(t)=0$ and $\ell_i(t)=r_i(t)$ remain frozen at their value at $\beta_i$. 
Convergence in the FSP definition is understood as convergence to $[\alpha_i, \beta_i]$ and Skorohod convergence
to $[\ell_i(t), r_i(t), \tau_i(t)], t\in [0,T],$ for each $i$, and the Skorohod convergence to $f(\cdot,t), ~t\in [0,T],$.

For each cluster $i$, we have the (same as in the case of A-perturbations) Lipschitz properties of $\ell_i(t), r_i(t), \tau_i(t)$ within $[\alpha_i,\beta_i]$. If we denote $\tau(t)=\sum_i \tau_i$, it is still Lipschitz.
We have exactly same properties describing the dynamics of each $\tau_i(t)$ in terms of the density at its left edge $\ell_i(t)$.

Finally, the following properties hold for $\rho < h$, and are proved analogously to the way it is done for A-perturbations: \\
(i) Denote by $U_i$ the set of those time points in $[0,\infty)$, where $\tau'_i(t)>0$. Then, 
$$
\sum_i \L(U_i) \le 1/(1+p). 
$$
(ii) For each cluster $i$, we have that $\alpha_i < \infty$ implies $\beta_i \wedge T - \alpha_i < 1$, i.e., the duration of any cluster is less than $1$.\\  
(iii) We have \eqn{eq-tau-int-bound}:
$$ 
\int_0^T \tau(t) dt \le \bar C = 1 + \frac{1}{2(1+p)}.
$$

\subsection{I-perturbations, low-density result.}
\label{sec-i-low-result}

\begin{theorem}
\label{th-I-low} 
Consider the model with I-perturbations and $\rho < h$.
Then, as $n\to\infty$, the sequence of the stationary distributions of
 $[f^n(\cdot,t), \tau^n(t)]$ is such that its any subsequential limit is concentrated on the 
 states with $\tau=0$ and finite number of clusters.
In other words, $\tau^n(\infty) \Rightarrow 0$ [equivalently, $\E \tau^n(\infty) \to 0$] and, consequently, 
the limiting flux (typical flux) $\phi(\rho) = \rho$. 
 \end{theorem}

{\em Proof of Theorem~\ref{th-I-low}} is essentially same as that of Theorem~\ref{th-A-low}. We do use the fact that both the pre-limit process and the FSPs have a finite number of clusters on any finite interval -- this means that, just like for A-perturbations, convergence of pre-limit trajectories to an FSP implies $\tau^n(\cdot)\to \tau(\cdot)$ u.o.c..
$\Box$

\subsection{I-perturbations, high density result and conjecture.}
\label{sec-i-pert}

It is natural to conjecture that for I-perturbations, $h< \rho < 1/2$, the same result as Theorem~\ref{th-A-high}
for A-perturbations holds.

\begin{conjecture}
\label{conj-high-dens}
Consider the model with I-perturbations and $h< \rho < 1/2$.
Then, 
$$
\phi(\rho) = \phi^*.
$$
 \end{conjecture}
 
 The proof of Theorem~\ref{th-A-high}, however does not completely go through. The key difficulty is that, for I-perturbations
 with $h< \rho < 1/2$, the number of clusters in steady-state is {\em not} uniformly stochastically bounded. The key ideas of the 
 proof of Theorem~\ref{th-A-high} can be applied to obtain the following partial result, which shows that, when $\rho > h$, 
 the limiting flux is strictly less than $\rho$.
 
 \begin{theorem}
\label{th-perturb}
For any $\epsilon>0$, there exists $\delta>0$, such that for all large $n$ and the density $\rho \in (h+\epsilon, 1/2)$, the flux $\phi(\rho,n) < \rho-\delta$.
 \end{theorem}
 
 We do not give a proof of this partial result, because it is a much more involved version of the proof of Theorem~\ref{th-A-high}, and does not provide new insights.

\begin{figure}
  \centering
  \begin{subfigure}[b]{0.4\linewidth}
    \includegraphics[width=\linewidth]{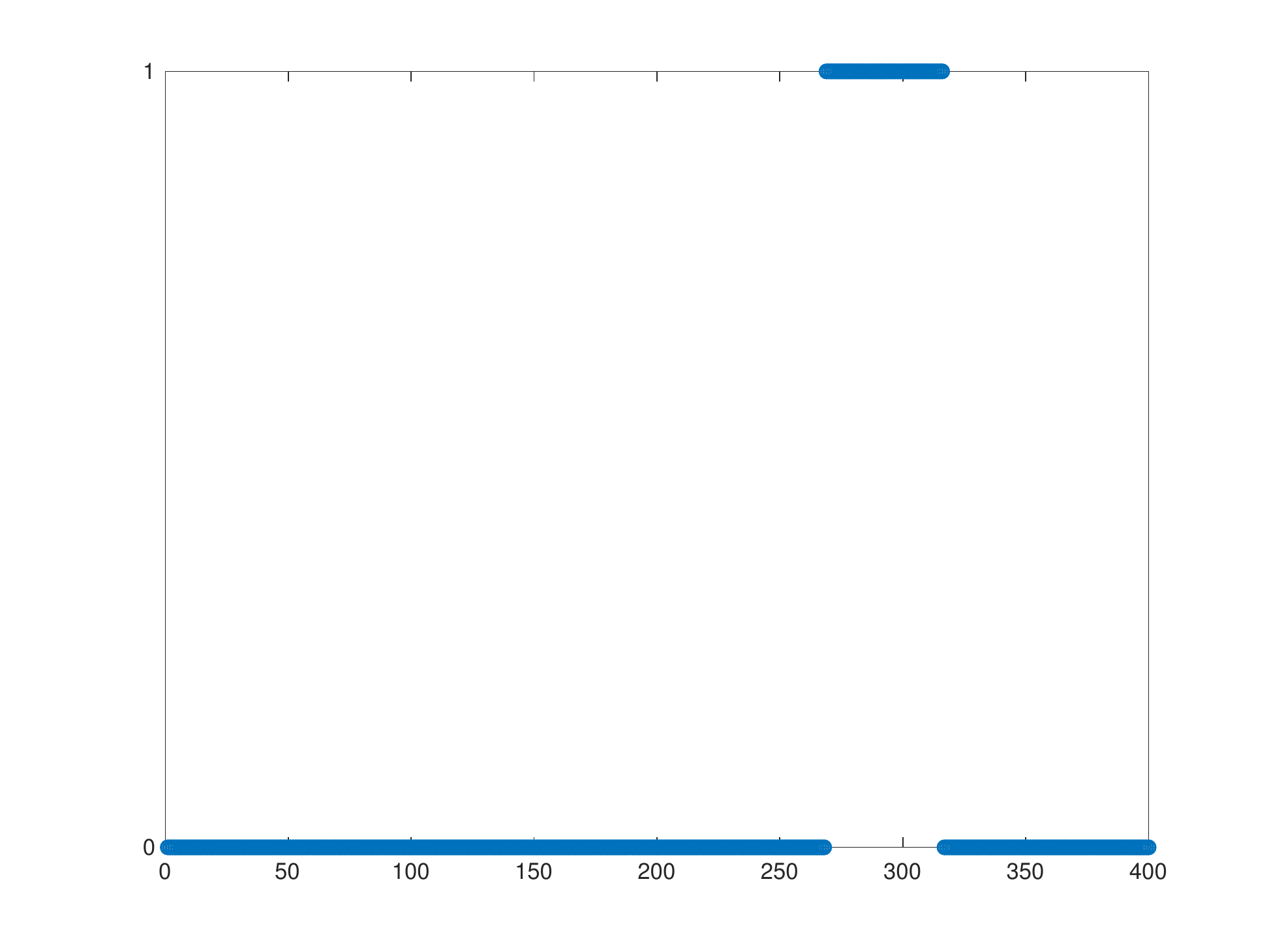}
  \end{subfigure}
  \begin{subfigure}[b]{0.4\linewidth}
    \includegraphics[width=\linewidth]{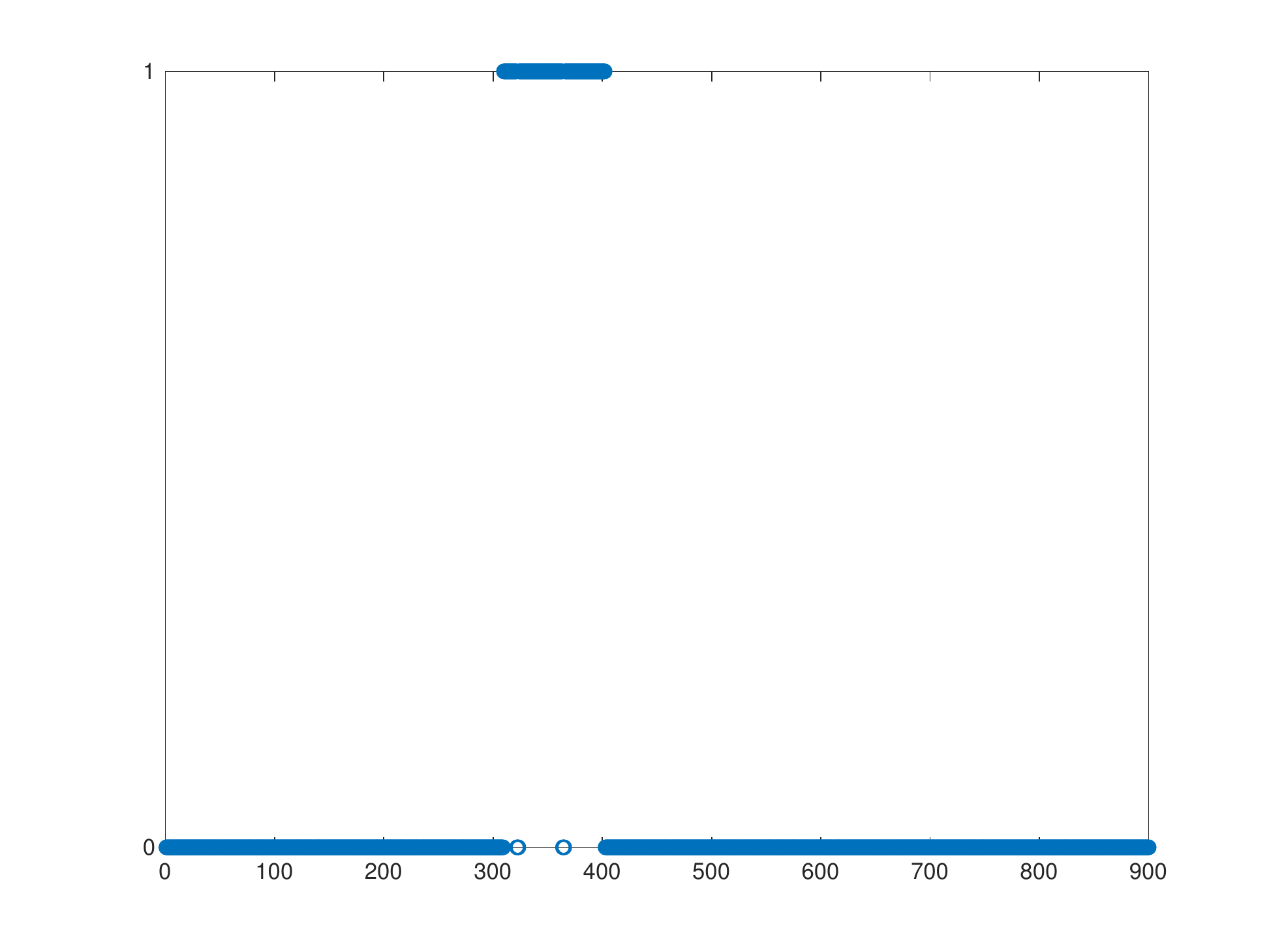}
  \end{subfigure}
  \caption{The states of the system with $n=400$ (left) and $n=900$ (right) sites, after $10^7$ time slots}
  \label{fig:tasep_h_i}
\end{figure}

Conjecture~\ref{conj-high-dens} is supported by extensive simulations. The simulations confirm that when when $\rho > h$ the limiting flux 
is indeed $\phi^*$. Furthermore, they show that the system ``typical'' state is ``close to'' MES in that it contain essentially a single cluster, ``punctured'' by holes due to perturbations. Figure \ref{fig:tasep_h_i} contains snapshots of the states of two systems (with $n=400$ sites in the left plot and $n=900$ sites in the right plot) with the values $\lambda = 1$, $p=1/2$ and $\rho = 0.4$ (so indeed $1/2 > \rho > h = p/(1+p)$). The systems are started from a random state, where the initially occupied sites are a random sample of $\rho n$ sites out of $n$ without replacement, and the snapshots are taken after $10^7$ time slots. For each site $i$ the plots show $x(i)*x(i+1)$, where $x(i)$ is $1$ or $0$ if
site $i$ is a particle or a hole, respectively. This is a convenient way to observe clusters in a large system. One can see that there is a single cluster on the left and essentially one cluster on the right, except it is ``punctured'' by two holes.

\section{A related zero-range model}
\label{sec-zrp}

Consider the following version (generalization) of TASEP-H. As before, there are $\rho n$ particles moving clockwise on $n$ sites forming a circle. There are no perturbations. At each time, the following occurs. A particle that has holes as both neighbors, moves right with probability $\pi\in (0,1]$. A particle that has another particle as the left-neighbor and a hole as the right-neighbor, moves right with probability $p\in (0,1)$. This model, as we will see shortly, can be viewed as a {\em zero-range model}, considered 
in \cite{Evans1997}. Denote by $v_{\pi}(\rho;n)$ and
$\phi_{\pi}(\rho;n)$
the steady-state particle velocity and flux, respectively, of this system. 
(Subscript $\pi$ will indicate that we refer to a flux for this model, as opposed to TASEP-H.)
It is shown in \cite{Evans1997} that the stationary distribution of this process has a product-form; moreover, the expressions are derived, which allow to determine the limiting steady-state particle velocity 
and limiting flux,  
$$
v_{\pi} = v_{\pi}(\rho) = \lim_{n\to\infty} v_{\pi}(\rho;n), ~~~ \phi_{\pi} = \phi_{\pi}(\rho) = \lim_{n\to\infty} \phi_{\pi}(\rho;n),
$$ 
in the same asymptotic regime as in this paper, with $n\to\infty$ and particle density $\rho$ staying constant. 
If $\pi=1$, we obtain exactly the basic TASEP-H model (without perturbations), so that $\phi_1 = \rho$ for any $\rho<1/2$.
If $\pi = 1-\epsilon < 1$, but close to $1$, then the zero-range model can be viewed as TASEP-H {\em with ``perturbations'' of a different kind}; namely,
at each time, with small probability $\epsilon>0$, a free particle stays in place instead of moving right. Note, however, that in the zero-range model (with $\pi$ independent of $n$) the ``perturbation'' rate is much higher -- $O(1)$ per particle per unit time -- than for our TASEP-H with perturbations model, where the perturbation rate is 
$O(1/n^2)$ per particle per unit time. 

Using the results of \cite{Evans1997}, we will prove the following

\begin{theorem}
\label{th-zrp} 
We have
\begin{equation} \label{eq:evans_eq1}
\lim_{\pi \to 1} \phi_{\pi}(\rho) = \begin{cases}
\rho, \quad \rho < h, \\
\phi^* = p(1-\rho), \quad \rho > h.
\end{cases}
\end{equation}
In other words, $\lim_{\pi \to 1} \phi_{\pi}(\rho) =\phi(\rho)$, where $\phi(\rho)$ is the limiting flux of TASEP-H with A-perturbations.
In particular, for $h< \rho < 1/2$, $\phi_{\pi}(\rho)$ is discontinuous at 
 $\pi=1$.
 \end{theorem}

{\em Proof. } 
Consider the equivalent model with holes moving left: a hole with exactly one consecutive particle immediately to the left, moves left with probability $\pi$; a hole with two or more consecutive particles immediately to the left, moves left with probability $p$. For convenience, denote by $\gamma = 1-\rho$ the density of holes.

Let $\eta_{\pi}=\eta_{\pi}(\rho)$ denote the limiting (in $n\to\infty$) steady-state velocity of a hole.
The relation to the velocity of a particle $v_{\pi}$ is via the equality of fluxes: $(1-\rho) \eta_{\pi} = \rho v_{\pi}$.

It is clear that \eqref{eq:evans_eq1} it equivalent to
$$
\lim_{\pi \to 1} \phi_{\pi} = \begin{cases}
1-\gamma, \quad \gamma > 1-h, \\
p \gamma, \quad \gamma < 1-h.
\end{cases}
$$
This, in turn, is equivalent to 
$$
\lim_{\pi \to 1} \eta_{\pi} = \begin{cases}
\frac{1-\gamma}{\gamma}, \quad \gamma > 1-h, \\
p, \quad \gamma < 1-h.
\end{cases}
$$
Let us use notation $q = \dfrac{1-\gamma}{\gamma}$. Note that as $h = p/(1+p)$, the condition $\gamma > 1-h$ is equivalent to $q < p$. So finally, we conclude that in order to show \eqref{eq:evans_eq1}, it is sufficient to show that
\begin{equation} \label{eq:evans_eq2}
\lim_{\pi \to 1} \eta_{\pi} = \begin{cases}
q, \quad q < p, \\
p, \quad q > p,
\end{cases}
\end{equation}
or simply $\lim_{\pi \to 1} \eta_{\pi} = q \wedge p$.

Results of \cite{Evans1997} are directly applicable to the model seen as movement of holes. In the notation of \cite[Section 5]{Evans1997}, $p(1)=\pi$, $p(k)=p$ for $k\ge 2$.
The function $f(k)$ (as in \cite[Section 5]{Evans1997}) in our case is then as follows (we can get rid of the common factor $1-\pi$, because this function is defined up to a positive factor):
$$
f(0) = 1,  ~~f(1) = \frac{1}{\pi},
$$
$$
f(k) = \frac{1-\pi}{\pi p} \left(\frac{1-p}{p}\right)^{k-2}, ~k\ge 2.
$$
Note that the generating function of $\{f(k)\}$ is
\beql{eq-genf}
G(z) = \sum_{k=0}^\infty f(k) z^k = 1 + \frac{z}{\pi} + \frac{1-\pi}{\pi p}  \frac{z^2}{1-(\frac{1-p}{p})z};
\eeql
clearly, for real $z\in [0,p/(1-p))$, $G(z)$ is increasing and $G(z)\uparrow \infty$, $z\uparrow p/(1-p)$,
and it is easy to check that same is true for the function $z G'(z)/G(z)$.

From \cite{Evans1997} the (limiting in $n \to \infty$) {\em fugacity} $z^*=z^*(\gamma,\pi)$ satisfies the equation
\beql{eq-fuga}
1-\gamma = \gamma z^* \frac{\partial \log G(z^*)}{\partial z^*} =  \gamma \frac{z^* G'(z^*)}{ G(z^*)}.
\eeql
Note that \eqref{eq-fuga} defines the fugacity $z^*$ uniquely, and $z^*$ is also clearly  continuous w.r.t. $\pi\in (0,1)$ and $(1-\gamma)/\gamma \in (0,\infty)$, viewed as parameters. (We will only use the continuity in $\pi$.)

The expression for the velocity $\eta_{\pi}$ in terms of fugacity $z^*$, is as follows (analogous to (11) 
in \cite{Evans1997}, holds for $z=z^*$):
$$
\eta_{\pi} = [\pi z f(1) + p \sum_{k=2}^\infty z^k f(k)] G(z)^{-1}.
$$
The equation above gives the dependence of $\eta_{\pi}$ on $z^*$. For the zero-range model that we consider, this dependence is very simple, and it holds
for a quite general function $p(k)$ (with $p(0)=0$ and $0<p(k)\le \delta < 1$ for $k\ge 1$). 
We have
$$
f(0) = 1,
$$
$$
f(k) = \frac{1}{1-p(k)}\prod_{m=1}^k \frac{1-p(m)}{p(m)}, \quad k \ge 1.
$$
Note that
$$
p(k) f(k) = \frac{p(k)}{1-p(k)}\prod_{m=1}^k \frac{1-p(m)}{p(m)} = \prod_{m=1}^{k-1} \frac{1-p(m)}{p(m)}  = (1-p(k-1))f(k-1).
$$
Hence the expression for $\eta_{\pi}$ is (we write $z$ instead of $z^*$ for simplicity below)
\begin{align*}
\eta_{\pi}  =  \frac{\sum_{k=1}^\infty p(k) f(k) z^k}{G(z)} = \frac{\sum_{k=1}^\infty (1-p(k-1)) f(k-1) z^k}{G(z)} \\ 
 =  \frac{z\sum_{k=1}^\infty f(k-1)z^{k-1} - z\sum_{k=1}^\infty f(k-1)p(k-1)z^{k-1}}{G(z)} \\ 
 =  \frac{zG(z) - z\eta_{\pi}G(z) }{G(z)} = z-z\eta_{\pi}.
\end{align*}
Therefore,
\beql{eq-vel}
\eta_{\pi} = \frac{z^*}{z^*+1}, ~~~\mbox{or} ~~~ z^* = \frac{\eta_{\pi}}{1-\eta_{\pi}}.
\eeql
Note that the dependence of $\eta_{\pi}$ on $\pi\in (0,1)$ (and $\gamma \in (0,1)$) is continuous.

Using \eqn{eq-fuga} and \eqn{eq-vel}, for a fixed $\gamma$,
we will now study the limit of $\eta_{\pi}$ as $\pi \uparrow 1$, and will prove  \eqref{eq:evans_eq2}.

Substituting \eqn{eq-genf} into \eqn{eq-fuga}, and then $z^*=\eta_{\pi}/(1-\eta_{\pi})$, we obtain 
\beql{eq-key-vh-p1}
(q-\eta_{\pi})(p-\eta_{\pi})(\eta_{\pi}(p-\pi-\pi p)+\pi p)= p^2(1-\pi)\eta_{\pi}(1-\eta_{\pi}).
\eeql
Note that, as $\pi\uparrow 1$, the RHS of \eqn{eq-key-vh-p1} converges to $0$, {\em while staying positive}, because for $\pi<1$ velocity cannot be $1$, and cannot be $0$ (in view of the LHS of \eqn{eq-key-vh-p1}).

Note that $\eta_{\pi}\le q$ must hold for any $\pi$, because $q$ is the max velocity achievable for the holes, achieved when all particles move at speed $1$. Therefore, to show that $\eta_{\pi} \le q \wedge p$, it suffices to show that 
$\eta_{\pi} \le p$ in the case $p<q$. If $\pi=p$, we have $\eta_{p} \le p$ simply because the velocity of a hole in this case cannot possibly exceed $p$. It is also easy to see that $\eta_{p}<p$. (In fact, in this case the velocity is known explicitly: $\eta_p = (1-\sqrt{1-4p\gamma(1-\gamma)})/(2\gamma)$, cf. \cite{Evans1997}.)
As we continuously increase $\pi$ in the interval $[p, 1)$, the RHS of \eqn{eq-key-vh-p1} must stay positive, so for all those $\pi$ we must have $\eta_{\pi}<p$.
This completes the proof of the fact that $\eta_{\pi} \le q \wedge p$ for all $\pi<1$.

Finally, again, as we continuously increase $\pi$ in the interval $[p, 1)$, velocity $\eta_{\pi}$ changes continuously,
the RHS of \eqn{eq-key-vh-p1} must stay positive and converge to $0$.
The only option is that $\eta_{\pi} \to q \wedge p$ as $\pi\uparrow 1$. $\Box$

\section{Discussion and further conjectures}
\label{sec-discuss-conj}

Our analysis of the TASEP-H model with perturbations suggests a number of generalizations and extensions. We discuss some of them.

\subsection{More general perturbation rates}

Consider TASEP-H with A-perturbations. Recall that for this process there are no perturbations as long as a cluster exists.
As $n\to\infty$, the average time it takes the cluster in a MES to dissolve, ``should'' grow exponentially fast with $n$. 
(A standard large-deviations analysis should apply.)
Therefore, if the inter-perturbation times grow sub-exponentally in $n$, then the fraction of time when the system spends in MES will dominate.
This basic intuition leads to the following Conjecture~\ref{conj-general}.

We will say that a positive non-increasing function $g(k) \to 0$, $k\to \infty$ is sub-exponential, or write $g(k) = \mbox{SUBEXP} (k)$, if
$$
\lim_{k\to\infty} \log(g(k))/k = 0.
$$

\begin{conjecture}
\label{conj-general}
Consider the model with either A-perturbations or I-perturbations.
Suppose, the probability of a perturbation (under the corresponding conditions for A- and I-perturbations) at a given time is
$\rho n g(n)$, where $g(n) = \mbox{SUBEXP}(n)$. ($g(n)$ is the perturbation rate, per particle per time unit.)
Then the limiting flux (typical flux) is as follows. \\
(i) If $0 < \rho < h$,
$$
\phi(\rho) = \lim_{n\to\infty} \phi(\rho;n) = \rho.
$$
(ii) If $h< \rho < 1/2$,
$$
\phi(\rho) = \lim_{n\to\infty} \phi(\rho;n) = \phi^*.
$$
\end{conjecture} 

The intuition that leads to Conjecture~\ref{conj-general} for TASEP-H with perturbations, as well as the informal connection between TASEP-H and the zero-range model, also suggests that the following conjecture is very plausible.

\begin{conjecture}
\label{conj-general-zrp}
Consider the zero-range model.
Suppose, parameter $\pi$ depends on $n$ so that $1-\pi(n) = \mbox{SUBEXP}(n)$.
Then the following holds. \\
(i) If $0< \rho < h$,
$$
\lim_{n\to\infty} \phi_{\pi(n)}(\rho;n) = \rho.
$$
(ii) If $h< \rho < 1/2$,
$$
\lim_{n\to\infty} \phi_{\pi(n)}(\rho;n)  = \phi^*.
$$
\end{conjecture}

\subsection{
Slow-to-start model}
\label{sec-sts}

Consider the following model, which is a discrete-time version of the {\em slow-to-start} (STS) model, studied in \cite{CFP2007}. 
The model is exactly as TASEP-H, except for the following. When a cluster of particles ``dissolves'' -- i.e., the cluster has consisted of exactly two particles, its ``right'' particle breaks away, and no particle joins the ``left'' particle from the left -- the ``left'' particle (which now has holes on both sides) does not become free immediately; instead, at each consecutive time it moves forward and becomes free with probability $p$. 

The STS model is primarily motivated by the movement of cars on a road. A car moves at a constant speed (normalized to $1$), unless it is blocked by a stopped car ahead; when the car is unblocked, it resumes movement not immediately, but after some delay (hence, ``slow start'').
Paper  \cite{CFP2007} studies the continuous STS model on the real line, with particles (cars) represented by points moving continuously at constant speed, unless stopped by running into a not moving particle ahead, and with exponentially distributed restart times. The results of \cite{CFP2007} concern the system evolution from an initial state, where none of the particles moves, towards a state where (in a sense) all particles move without bumping into each other. (Loosely speaking, the latter state is an absorbing state in the terminology of this paper.) 

Considering the typical flux of the discrete-time STS model, described above, is of interest for the same reasons as for TASEP-H. The STS model with A-perturbations is defined the same way as for TASEP-H, i.e. a perturbation is applied with certain probability when all particles are free. (Recall that for STS model a particle is free when it is not only has holes on both sides, but also is actually moving.) The STS model with I-perturbations is defined exactly as for TASEP-H. 

It is not hard to check that the main results of the present paper for TASEP-H with A- and I-perturbations hold {\em as is} for the discrete-time STS model with A- and I-perturbations. Same proofs apply, with slight modifications.

\appendix 
\section{Proof of Proposition~\ref{th-ballot}.}
\label{app-ballot}

The result is derived  from Theorem 1.2.5 -- a ballot theorem -- in \cite{Tac1967}.

Without loss of generality, we can assume that $m=n$. (Otherwise, $m$ and all $k_r$ can be rescaled by the same factor, and this obviously does not change  $N$.) Then, conditions \eqn{eq-bal1} and \eqn{eq-bal2} become, respectively,
\beql{eq-bal111}
\psi(j+r) -\psi(j) <  r, ~~~  r=1,\ldots,n,
\eeql
and
\beql{eq-bal222}
N \ge \lceil n- \psi(n)   \rceil.
\eeql
For each integer $j\ge 0$, let $\eta(j) = 1$ if \eqn{eq-bal111} holds, and $\eta(j) = 0$ otherwise. (So that, $N = \sum_0^{n-1} \eta(j)$.)
Let us extend the definition of $\psi(j)$ to all real $u\ge 0$ by letting $\psi(u) = \psi(\lfloor u \rfloor)$. For each $u\ge 0$, let $\delta(u)=1$ if
\beql{eq-delta-def}
v-u \ge \psi(v)-\psi(u) ~~\mbox{for all $v\ge u$},
\eeql
and $\delta(u)=0$ otherwise. Denote by $D$ the set of those $u \in [0,n)$, for which $\delta(u)=1$. By Theorem 1.2.5 in \cite{Tac1967}, 
$\L(D) = n-\psi(n)$, where $\L$ denotes the Lebesgue measure. For any interval $[j, j+1)$, $j=0,1,\ldots,n-1$, 
observe the following: if $\delta(u)=1$ for some $j < u <j+1$, then necessarily $\eta(j)=1$. 
Indeed, in this case, for any $r=1,\ldots,n$,
$$
r- [\psi(j+r) -\psi(j)] = (u-j) + (j+r-u) -  [\psi(j+r) -\psi(u)] \ge u-j >0,
$$
where the first inequality is by \eqn{eq-delta-def}.
This implies that for each interval $[j, j+1)$ such that $\L(D \cap [j, j+1))>0$ we have $\eta(j)=1$. Since, obviously, $\L(D \cap [j, j+1))\le 1$ for any $j$, and 
$$
\sum_{j=0}^{n-1} \L(D \cap [j, j+1)) = \L(D) = n-\psi(n),
$$
we see that the (integer) number $N$ of those $j \in \{0,1,...,n-1\}$ for which $\eta(j)=1$ cannot be less than $n-\psi(n)$. 
This proves \eqn{eq-bal2}.
$\Box$

\section{A CSMA model indirectly motivating the holdback property of TASEP-H}
\label{subsec:stork}

In this section we describe an interacting particle system, which is a model of a fairly realistic wireless network under a CSMA protocol. 
This system is far more complicated than TASEP-H, but exhibits qualitatively similar behavior. We will formulate some hypotheses and discuss the difficulties one faces in analyzing such a system. We briefly introduce the model here and refer to \cite{ShSt2018} for an account of known results on the topic, in particular stability analysis. 

As in the TASEP-H setting, there are $n$ sites arranged in a circle, and $\rho n$ particles, moving clockwise (``right'') in 
discrete time. 
(A time index is sometimes referred to as a ``time slot.'') 
The particle density $\rho\in (0,\infty)$ (not necessarily less than $1$).
At each time a particle may move to its right-neighbor site.
There is no restriction on the number of particles present at a site at a given time. A particle may move to the right-neighbor site even if that site has other particles. The particles that do move at a given time are chosen according to a randomized competition described as follows. At  each time the sites are given access priorities, forming a permutation of numbers $1,\ldots,n$,
picked independently (across time slots), uniformly at random from all possible permutations. From the site with the highest priority, if it is not empty (has one or more particles), exactly one particle (chosen at random, say) moves to its right-neighbor site.  From the site with the second-highest priority, exactly one particle moves to the right-neighbor site, as long as the site is not empty and that no particle has already moved from any of its two neighbor sites. And so on until all sites are checked in their priority order. The procedure described above is repeated independently over time slots.

The restriction that a particle can only move from a site if no particle moved from any of the neighbors of the site models the interference in wireless networks, and it is the most important characteristic of models of these networks.

\begin{figure}[h!]
\centering
\includegraphics[width=0.7\linewidth]{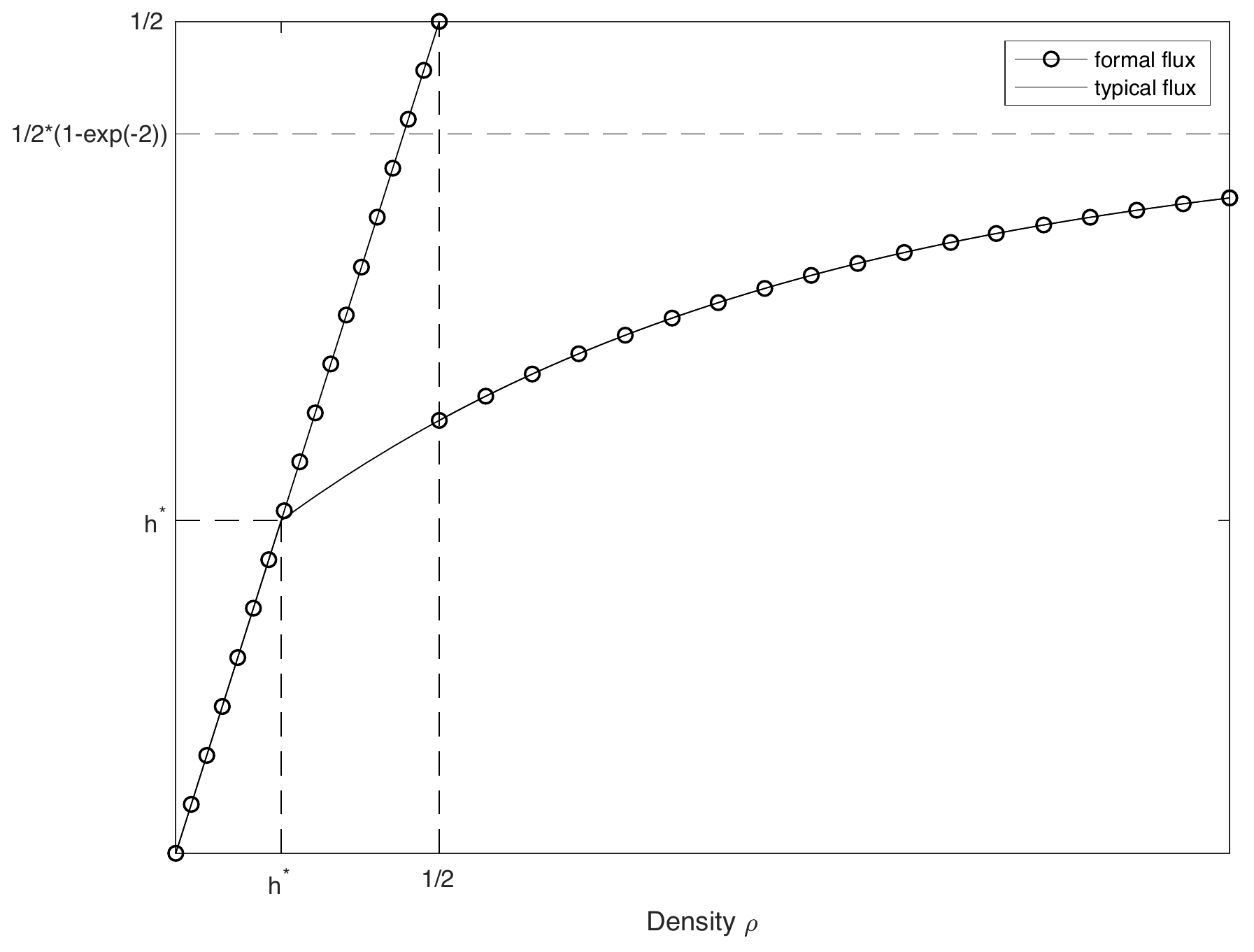}
\caption{Formal and typical flux in the CSMA model, versus density $\rho$, as suggested by simulations. For $\rho< h^*$ and $\rho>1/2$ the formal and typical fluxes coincide}
\label{fig:flux_stork}
\end{figure}

We are interested in the flux of this system. As in the TASEP-H model, if $\rho < 1/2$ an ideal (absorbing) state will be eventually reached and the system will stay in it, thus the formal flux is equal to $\rho$ as long as $\rho < 1/2$ (see Figure~\ref{fig:flux_stork}). 

A typical flux is defined analogously to that for TASEP-H, i.e. as the limit of the flux of the system with perturbations, as $n\to\infty$.
The typical flux that we observe in simulations is shown in Figure~\ref{fig:flux_stork}. 
(The typical flux plot in Figure~\ref{fig:flux_stork} is just a qualitative illustration -- it does not describe a specific simulation experiment.)
We observe that the typical flux is monotone increasing, asymptotically converging to the {\em parking constant} $c=(1/2)(1-e^{-2}) \approx 0.43$ as $\rho \to \infty$. The asymptotic limit $c$ is intuitive. Indeed, as $\rho \to \infty$, almost all sites will be occupied most of the time. If so, 
when $n$ is large, the fraction of sites that will move a particle at a given time will be equal to the expected number of cars parked per space-slot in the so-called discrete parking system (we refer to \cite{ShSt2018} for an explanation of the connections between the models).
We also {\em conjecture} that there is a phase transition at some density level $h^*$.  The conjectured phase transition should be qualitatively similar to that 
for TASEP-H model. Namely, when the density $\rho<h^*$, the clusters that may emerge will not have a tendency to grow very large, because, if/when they become large they  ``lose'' particles on the right faster than ``acquire'' new particle on the left (due to free particles joining it). As a result, when $\rho<h^*$ almost all (in the $n\to\infty$ limit -- all) particles remain free. If $\rho>h^*$, very large clusters form and contain a non-zero fraction of particles, while the density of particles in the sparse intervals is $h^*$. (The value of $h^*$ is not known to us. A simple heuristic argument leads us to believe that $h^* \ge 1/9$. Indeed, the rate at which particles ``break away'' on the right  from a cluster, followed by a sparse interval, is at least $1/8$; within any two time slots, a particle breaks away with probability at least $(1/2)^2$. And if the density of the sparse interval on the left of the cluster is less than $1/9$, then particles join the cluster on the left at the rate at most $(1/9)/(1-1/9)=1/8$. Thus, if the overall particle density $\rho< 1/9$, clusters should have a tendency to dissolve.)

We further conjecture that for $\rho > 1/2$ the formal flux coincides with the typical flux. That is,
the formal flux has a negative ``jump'' at $\rho=1/2$, and then is equal to the typical flux for $\rho > 1/2$ (see Figure~\ref{fig:flux_stork}).

We note that an analysis of the system in this section is much more difficult than that of the TASEP-H system due to a number of reasons, including, but not limited to, the following: at any time, any cluster may in principle break into any number of clusters even without perturbations; similarly, any number of clusters may merge and form a much larger cluster; even if the right-most site of a cluster moves a particle to the right, the latter particle does not necessarily break away from the cluster, because the site may have had other particles in it.

\end{document}